\documentclass[a4paper,11pt]{amsart}
\usepackage{graphicx}
\usepackage{caption}
\usepackage[latin1]{inputenc}
\usepackage[left=0.6in, top=0.8in, paperwidth=8in, paperheight=11in]{geometry}

\newtheorem{thm}{Theorem}[section]
\newtheorem{cor}[thm]{Corollary}
\newtheorem{lem}[thm]{Lemma}
\newtheorem{prop}[thm]{Proposition}
\theoremstyle{definition}
\newtheorem{defn}[thm]{Definition}
\theoremstyle{remark}

\numberwithin{equation}{section}

\begin{document}
\title[Nonrigidity for circle homeomorphisms with several break points
]{Nonrigidity for circle homeomorphisms \\ with several break points
 }
\author{Abdelhamid Adouani and Habib Marzougui}
\address{University of Carthage, Faculty of Science of Bizerte, Department of Mathematics, Jarzouna, 7021, Tunisia}
\email{arbi.abdelhamid@gmail.com}
\email{hmarzoug@ictp.it, habib.marzougui@fsb.rnu.tn}

\begin{abstract} Let $f$ and $g$ be two class $P$-homeomorphisms of the circle $S^{1}$ with break points singularities, that are
differentiable maps except at some singular points where the
derivative has a jump. Assume that the derivatives $\textrm{Df}$ and $\textrm{Dg}$ are
absolutely continuous on every continuity interval of $\textrm{Df}$
and $\textrm{Dg}$ respectively. Denote by $C(f)$ the set of break points of $f$. For $c\in S^{1}$, denote by
$\pi_{s, O_{f}(c)}(f)$ the product of $f-$ jumps in break points lying to the $f-$ orbit
of $c$ and by $\textrm{SO}(f) = \{O_{f}(c):~c \in C(f)~\textrm{and}~\pi_{s, O_{f}(c)}(f)\neq 1\}$,
called the set of singular $f$-orbits. The maps $f$ and $g$ are called break-equivalent if there exists a topological
conjugating $h$ such that $h(\textrm{SO}(f))=\textrm{SO}(g)~~ \textrm{and} ~~ \pi_{s,
O_{g}(h(c))}(g) = \pi_{s, O_{f}(c)}(f) ~~ \textrm{for all}~~ c\in
\textrm{SO}(f)$. Assume that $f$ and $g$ have the same
irrational rotation number of bounded type. We prove that if
$f$ and $g$ are not break-equivalent, then any topological conjugating $h$ between $f$
and $g$ is a singular function i.e. it is a continuous on $S^{1}$,
but $\textrm{Dh}(x)=0$ a.e. with respect to the Lebesgue measure. As a consequence if for some point $d\in \textrm{SO}(f)$, $\pi_{s, O_{g}(d)}(g)\notin \{\pi_{s, O_{f}(c)}(f):
c\in C(f)\}$, then the homeomorphism conjugation $h$ is a singular function. This later result generalizes
previous results for one and two break points obtained by
Dzhalilov-Akin-Temir and Akhadkulov-Dzhalilov-Noorani. Moreover, if $f$ and $g$ do not have the same number of singular orbits
then the homeomorphism conjugating $f$ to $g$ is a singular function.
We also deduce that under the case of rotation number of bounded type that
if $f$ does not have the (D)-property (i.e. $\pi_{s, O_{f}(c)}(f)= 1$, for all $c\in C(f)$) and $g$ has the
(D)-property, then the map conjugating $f$ to $g$ is a singular function. In particular, 
if $f$ does not have the (D)-property with rotation number of bounded type, then the invariant measure $\mu_{f}$
is singular with respect to the Lebesgue measure. Theses results cannot be extended to rotation number not of bounded type 
even for piecewise linear homeomorphisms (see very recently, an example due to Teplinsky in \cite{aT15}).

\end{abstract} \maketitle


\section{\bf Introduction}
\bigskip

Denote by $S^{1} = \mathbb{R}/\mathbb{Z}$ the circle and  $p :
\mathbb{R}\longrightarrow S^{1}$ the canonical projection. Let  $f$
 be an orientation preserving homeomorphism of $S^{1}$. The homeomorphism  $f$ admits a lift $\widehat{f} :
\mathbb{R}\longrightarrow \mathbb{R}$  that is an increasing
homeomorphism of  $\mathbb{R}$  such that $p\circ\widehat{f} =
f\circ p$. Conversely, the projection of such a homeomorphism of
$\mathbb{R}$ is an orientation preserving homeomorphism of
$S^{1}$. Let  $x\in S^{1}$. We call \emph{orbit} of  $x$  by  $f$
the subset $O_{f}(x) = \{f^{n}(x): n\in\mathbb{Z} \}$. Historically, the study of the dynamics of 
circle homeomorphisms was
initiated by Poincar\'e (\cite{hP}, 1886), who introduced the
rotation number of a homeomorphism $f$  of $S^{1}$ as $\rho (f) =
\underset{n\to +\infty}{\lim}\frac{\widehat{f}^{n}(\widehat{x}) -
\widehat{x}}{n}~(\textrm{mod } 1)$, where $\widehat{x}\in
\mathbb{R}$ such that $p(\widehat{x})= x$. Poincar\'e shows that
this limit exists and does not depend on neither $x$ nor the lift
$\widehat{f}$ of $f$. We say that $f$ is semi-conjugate to the
rotation  $R_{\rho(f)}$ if there exists an orientation preserving
surjective continuous map $h: S^{1}\longrightarrow S^{1}$ of degree
one such that $h\circ f = R_{\rho(f)}\circ h$.

{\it Poincar\'e's theorem.} Let  $f$  be an homeomorphism of $S^{1}$
with irrational rotation number $\rho(f)$. Then $f$ is
semi-conjugate to the rotation  $R_{\rho(f)}$.

A natural question is whether the semi-conjugation $h$ could be
improved to be a conjugation, that is $h$ to be an homeomorphism. In
this case, we say that $f$ is topologically conjugate to the
rotation $R_{\rho(f)}$. In this direction, Denjoy (\cite{aD32})
proves the following:
\medskip

{\it Denjoy's theorem \cite{aD32}}. Every $C^{2}$-diffeomorphism $f$
of $S^{1}$ with irrational rotation number $\rho(f)$ is
topologically conjugate to the rotation $R_{\rho(f)}$.
\medskip

Denjoy asked whether or not $C^{2}$-diffeomorphisms $f$ of $S^{1}$
are ergodic with respect to the Lebesgue measure $m$ ($f$ is said to
be ergodic with respect to $m$ if any $f$-invariant measurable set
$A$ has measure $m(A)$ equal to $0$ or $1$). Simultaneously, Herman
and Katok gave a positive answer to this question:
\medskip

{\it Herman-Katok's theorem} (\cite{yKdO89}, \cite{bHaK95}). Every
$C^{2}$-diffeomorphism  $f$  with irrational rotation number is
ergodic with respect to the Lebesgue measure $m$.
\medskip

It is well known that an homeomorphism  $f$ of  $S^{1}$  with
irrational rotation number preserves a unique normalized measure on
$S^{1}$, denoted by  $\mu_{f}$. If  $\widehat{h}$ is the lift of
$h$,  by taking  $\widehat{h}(0) = 0$,  the conjugating
homeomorphism  $h$ is unique and related to $\mu_{f}$ by
$h(x)=p(\mu_{f}([0,x])) \in S^{1}$ for $x \in S^{1}$. Uniqueness of
$\mu_{f}$ implies that $\mu_{f}$ is either singular, or absolutely
continuous with respect to  $m$, in this second case, $h$ is an
absolutely continuous function. Recall that $\mu_{f}$ is said to be
\textit{singular} with respect to the Lebesgue measure $m$ on
$S^{1}$ if there exists a measurable subset $E$ of $S^{1}$ such that
$\mu_{f}(E) = 1$ and $m(E) = 0$.
 In fact, if $\mu_{f}$  is absolutely continuous with respect to
 the Lebesgue measure $m$ and  $f$ is a  $C^{2}$-diffeomorphism,
 $\mu_{f}$ is necessarily equivalent to  $m$  as a consequence of Herman-Katok's theorem
 above (i.e. $m$ is absolutely continuous with respect to $\mu_{f}$
 and conversely).
\medskip
In the sequel we deote by $\mathbb{R}^{\ast}=\mathbb{R}\backslash
\{0\}$ and  $\mathbb{N}^{\ast}=\mathbb{N}\backslash \{0\}$.

\textit{Definition}. A real number  $\alpha\in ]0, 1[$ is called
\textit{Diophantine with exponent $\delta\geq 0$} if there is a
constant $c(\alpha)>0$ such that
$$ \ (1) \qquad \qquad \arrowvert\alpha -\dfrac{p}{q}\arrowvert\geq \dfrac{c(\alpha)}{q^{2+\delta}}\qquad \textrm{ for any }\  \dfrac{p}{q}\in \mathbb{Q}.$$
\medskip
\
\\
A number that is neither rational nor Diophantine is called a
\textit{Liouville number}. \
\\
 Every real number $\alpha\in ]0, 1[$ has a continued fraction expansion
represented by
$$\alpha = \frac{1}{a_{1}+ \frac{1}{a_{2}+ \dots}}: =
[a_{1},a_{2},\dots,a_{n},\dots]$$ \
\\
where  $a_{m} \in \mathbb{N}^{\ast}$, $m\in \mathbb{N}^{*}$ are
called \textit{partial quotients} of $\alpha$. When ($a_{m})_{m\in
\mathbb{N}}$ is bounded, $\alpha$ is said to be of \textit{bounded
type}. This is equivalent to the fact that $(1)$ holds with $\delta
= 0$.
\medskip

The problem of smoothness of the conjugacy $h$ of smooth
diffeomorphisms to rotations is now very well understood (see for
instance \cite{mH79}, \cite{yKdO89}, \cite{kmKyS89}, \cite{KT09},
\cite{jcY84}). 



\medskip



We refer the reader to the books \cite{aN10} and \cite{wDvS91} for a
thorough account on circle homeomorphisms.
\smallskip

The situation is more complicated for circle homeomorphisms with
break points or shortly, class $P$-homeomorphisms (see the
definition below). This class are known to satisfy the conclusion of
Denjoy's theorem (Corollary 2.6) (see also \cite{yKdO89};
\cite{mH79}, chapter VI) and, with additional regularity, of
Herman-Katok's theorem (see \cite{kaA06}). However,
Katznelson-Ornstein's theorem \cite{yKdO89} cannot be extended in
general to class $P$. The study of the regularity of the invariant
measures of class $P$-homeomorphisms arises then naturally.
\bigskip

\textbf{Class  $P$-homeomorphisms} \\
The following definition is du to M.R. Herman.
\begin{defn}(see \cite{mH79}, p.74) An orientation preserving homeomorphism $f$ of $S^{1}$ is called a \emph{class  $P$}-homeomorphism if it is differentiable
except at finitely many points, the so called \textit{break points}
of $f$, at which left and right derivatives (denoted, respectively,
by $\textrm{Df}_{-}$ and $\textrm{Df}_{+}$) exist and such that the
derivative  $\textrm{Df}: S^{1} \longrightarrow \mathbb{R}_{+}^{*}$
has the following properties:
\begin{itemize}
\item [\textbullet] there exist two constants  $0<a<b<+\infty$  such that: $a<\textrm{Df}(x)<b,$  for every $x$, where $\textrm{Df}$  exists,

\item [\textbullet]  $a <\textrm{Df}_{+}(c)< b$ and  $a<\textrm{Df}_{-}(c)< b$  at the break points $c$.

\item [\textbullet]  $\log \textrm{Df}$  has bounded variation on  $S^{1}$.
\end{itemize}
\end{defn}
\medskip

We pointed out that the third condition implies the two ones.
\medskip
The ratio $\sigma_{f}(c): =
\dfrac{\textrm{Df}_{-}(c)}{\textrm{Df}_{+}(c)}$  is called the
\emph{$f$-jump} in $c$.
\medskip

 Denote by
\begin{itemize}
 \item [\textbullet] $C(f) = \{c_{0}, c_{1}, c_{2},\dots, c_{p}\}$  the set of break points of $f$ in
$S^{1}$

 \item [\textbullet] $c_{p+1}: = c_{0}$.

 \item [\textbullet] $\pi_{s, O_{f}(c)}(f): =  \underset{x\in C(f)\cap O_{f}(c)}{\prod}\sigma_{f}(x)$,
for every $c\in S^{1}$. \
\\
 \item [\textbullet] $\pi_{s}(f)$ the product of $f$-jumps in the break points of $f$:
$$\pi_{s}(f) =  \underset{c\in C(f)}{\prod}\sigma_{f}(c).$$
\\
\item [\textbullet] $SO(f)=\{O_{f}(c):~c \in C(f)~\textrm{and}~\pi_{s, O_{f}(c)}(f)\neq 1\}$

\end{itemize}
\bigskip

The total variation of $\log \textrm{Df}$ we denote by $V =
\textrm{Var}(\log \textrm{Df})$. We have $$~V : \ = \
\sum_{j=0}^{p}\textrm{Var}_{[c_{j},c_{j+1}]}\log \textrm{Df} +
|\log(\sigma_{f}(c_{j}))|< +\infty.$$ In this case,  $V$ is the
total variation of  $\log \textrm{Df}$,  $\log \textrm{Df}_{-}$,
$\log \textrm{Df}_{+}$.
\bigskip

We notice the following properties:
\begin{itemize}
 \item [\textbullet] If $f$ is a class $P$-homeomorphism of $S^{1}$ which is $C^{1}$ on $S^{1}$ then $f$ is a $C^{1}$-diffeomorphism.

\item [\textbullet] If $f$ is a class $P$-homeomorphism of $S^{1}$ then $D^{2}f \in L^{1}(S^{1})$: Indeed, this follows from
(\cite{HS}, Theorem $18.14$, p. 284) since $\textrm{Df}$ has bounded
variation on $S^{1}$.

\item [\textbullet] If $f$ is a class $P$-homeomorphism of $S^{1}$ and $\textrm{Df}$ is absolutely continuous on every connected interval of $S^{1}\backslash C(f)$ then $f\in C^{1}(S^{1}
\backslash C(f))$.
\end{itemize}

 Among the simplest examples of class $P$-homeomorphisms, there are:
\medskip

\textbullet \; $C^{2}$-diffeomorphisms,

\textbullet \; Piecewise linear $\textrm{PL}$-homeomorphisms.

An orientation preserving homeomorphism $f$ of $S^{1}$ is called a
\emph{PL-homeomorphism} if $f$ is differentiable except at finitely
many break points $(c_{i})_{1\leq i\leq p}$ of $S^{1}$  such that
the derivative  $\textrm{Df}$ is constant on each  $]c_{i}, \
c_{i+1}[$.
\bigskip

Denote by
\begin{itemize}

\item [\textbullet] $\mathcal{P}(S^{1})$ the set of class $P$-homeomorphisms of $S^{1}$.

\item [\textbullet] $\textrm{PL}(S^{1})$ the set of orientation preserving piecewise linear $\textrm{PL}$-homeomorphisms of $S^{1}$.
\end{itemize}
\medskip
We notice that $\mathcal{P}(S^{1})$ is a group which contains
$\textrm{PL}(S^{1})$ (cf. \cite{am08}).
\bigskip

\begin{defn}\cite{am15}  Let  $f\in \mathcal{P}(S^{1})$. We say that $f$ has the \textit{($D$)-property}
if the product of $f$-jumps at the break points of each orbit is
trivial, that is  $\pi_{s, O_{f}(c)}(f)= 1$, for every $c\in C(f)$.
\end{defn}
\medskip

In particular, if $f$ has the ($D$)-property, then $\pi_{s}(f) = 1$.
Conversely, if all break points belong to the same orbit and
$\pi_{s}(f) = 1$ then  $f$ has the ($D$)-property.
 \bigskip

 Let  $f\in \mathcal{P}(S^{1})$. We say that $f$ satisfies the \textit{Katznelson-Ornstein (KO)} if
 $\textrm{Df}$ is absolutely continuous on every continuity interval of $\textrm{Df}$.


\medskip


%
%
\medskip

%
%

\textbullet ~~ Recently, there has been a significant progress in
the problem of the regularity of the conjugating map between two
class $P$-homeomorphisms with the \textit{same}
irrational rotation number. The case of two class $P$-homeomorphisms with one break point and with the
same jump, was studied by Khanin and
Khmelev in \cite{KhKm} and by Teplinskii and Khanin in \cite{TK}.

%
\medskip

\textbullet ~~ The case of two class $P$-homeomorphisms with \textit{one} break point and with
\textit{distinct} jumps, was studied by Dzhalilov, Akin and Temir in
\cite{AHS10}. Their results show that the conjugation map between two such maps is singular.

\textbullet ~~ The case of two class $P$-homeomorphisms with \textit{two} break points was studied by Akhadkulov,
Dzhalilov and Mayer in \cite{HAD11} (see Corollary \ref{c:14}).
\medskip

\textbullet ~~ The case of \textit{several} break point with
\textit{distinct total jumps}, was studied by Dzhalilov, Mayer and Safarov in
\cite{DMS}, and independently by the first author in \cite{A}.
Their results show that the conjugation map between two such maps is singular:

{\it Theorem $($\cite{A}, \cite{DMS}$)$} \textit{Let $f$ and $g$ be two class $P$-homeomorphisms of the circle $S^{1}$ and same
irrational rotation number. Assume that $f$ and $g$ satisfy that satisfy the Katznelson-Ornstein (KO).
If $\pi_{s}(f)\neq \pi_{s}(g)$ then the homeomorphism map $h$ conjugating $f$ and $g$ is a singular function i.e. it is continuous on $S^{1}$ and $\textrm{Dh}(x)=0$ a.e. with respect to the Lebesgue measure.}
\bigskip

\begin{defn}
We say that two class $P$-homeomorphisms $f$ and $g$ of the circle
$S^{1}$ are \textit{break-equivalent} if there exists a topological
conjugating $h$ such that
\begin{itemize}
 \item[(1)] ~~ $h(\mathrm{SO}(f)) = \mathrm{SO}(g)$ \ and
 \item[(2)] ~~ $\pi_{s, O_{g}(h(c))}(g) = \pi_{s, O_{f}(c)}(f)$, ~~ \textrm{for all}~~ $c\in
C(f)$
\end{itemize}
\end{defn}
\medskip

In particular, if all points of $C(f)$ (resp. $C(g)$) are on pairwise
distinct orbits then

- $\mathrm{card}(\mathrm{SO}(f)) = \mathrm{card}(C(f))$
(resp. $\mathrm{card}(\mathrm{SO}(g))=\mathrm{card}(C(g))$),

- $f$ and $g$ are break-equivalent if there exists a topological
conjugation $h$ from $f$ to $g$ such that:
\begin{itemize}
 \item[(1)] ~~ $h(C(f))=C(g)$ \ and
 \item[(2)] ~~ $\sigma_{g}(h(c))=\sigma_{f}(c) ~~ \textrm{for all} ~~ c \in C(f).$
\end{itemize}
\bigskip

This definition coincides with that of \cite{KcDs}. It is easy to see
that if there is a $C^{1}-$ conjugation between $f$ and $g$ then $f$
and $g$ are break-equivalent. The case of break-equivalent
$C^{2+\alpha}$ homeomorphisms $f$ and $g$ with
$\pi_{s}(f)=\pi_{s}(g)=1$ and some combinatorial conditions was
study by Cunha-Smania in \cite{KcDs}. It was proved that any two
such homeomorphisms are $C^{1}-$ conjugated.
\bigskip
\\
In the present paper we consider \textit{non break equivalent}
class $P$-homeomorphisms having \textit{several break points} with
coinciding irrational rotation number of bounded type. The main purpose is to prove the following:

\
\\
{\bf Main Theorem }(Nonrigidity). \textit{Let $f$ and $g$ be two
class $P$-homeomorphisms of the circle $S^{1}$ with several break
points that satisfy the Katznelson-Ornstein (KO) and they have the
same irrational rotation number of bounded type. If $f$ and $g$ are not break equivalent then
the homeomorphism $h$ conjugating $f$ to $g$ is a singular function
i.e. it is continuous on $S^{1}$ and $\textrm{Dh}(x)=0$ a.e. with
respect to the Lebesgue measure}.
\medskip

As a consequence we have
\medskip

\begin{cor}\label{c:12}  Let $f$ and $g$ satisfy the assumptions of the main theorem. If \textrm{card} $(SO(f))\neq $ \textrm{card} $(SO(g))$ then
the map conjugating $f$ to $g$ is a singular function. In particular, this hold if $f$ (resp. $g$) has exactly $n$
break points (resp. $m$ break points) lying on pairwise distinct 
$f$-orbits (resp. $g$-orbits) with $n, m\in \mathbb{N}$, $n\neq m$.
\end{cor}
\bigskip

\begin{cor}\label{c:13}  Let $f$ and $g$ satisfy the assumptions of the main theorem. If there is some point $d \in SO(f)$
such that $\pi_{s,O_{g}(d)}(g) \notin \{\pi_{s,O_{f(c)}}(f): ~~c \in
C(f)\}$ then the map $h$ conjugating $f$ to $g$ is singular.
\end{cor}
In particular, for \textit{non break equivalent} homeomorphisms with
two break points, we get, as consequence:
\medskip

\begin{cor}\label{c:14} $($Akhadkulov-Dzhalilov-Noorani's theorem$)$ \cite{ADN}
Let $f_{i} \in \mathcal{P}(S^{1}),~i=1,2$ be circle homeomorphisms with
two break points $a_{i},~b_{i}$. Assume that $f_{i},~i=1,2$ satisfy
the Katznelson-Ornstein (KO) condition and $D^{2}f_{i}\in
L^{p}(S^{1})$ for some $p>1$. Assume that:
\begin{itemize}
  \item [(i)] The rotation numbers $\rho(f_{i})$ of $f_{i},~i=1,2$
are irrational of bounded type and coincide
$\rho(f_{1})=\rho(f_{2})=\rho, ~\rho \in \mathbb{R}\backslash
\mathbb{Q}$,
  \item [(ii)] $\sigma_{f_{1}}(a_{1}) \notin
\{\sigma_{f_{2}}(a_{2}),~\sigma_{f_{2}}(b_{2})\}$,
  \item [(iii)]
$\sigma_{f_{1}}(a_{1})~\sigma_{f_{1}}(b_{1})=\sigma_{f_{1}}(a_{2})~\sigma_{f_{1}}(b_{2})$,
  \item [(iv)] The break points of $f_{i},~i=1,2$ do not lie on the
same orbit,
\end{itemize}
Then the map $h$ conjugating $f_{1}$ to $f_{2}$ is singular.
\end{cor}
\medskip

\begin{cor}\label{c:15}
Let $f$ and $g$ satisfy the assumptions of the main theorem. Then:
\begin{itemize}
    \item [(i)] If $g$ does not have the $(D)-$ property and $f$ has
the $(D)-$ property then the map conjugating $f$ to $g$ is
singular.
    \item [(ii)] If $f$ and $g$ have the $(D)-$property and $D^{2}f,~D^{2}g \in
L^{p}(S^{1})$ for some $p>1$ then the map conjugating $f$ to $g$ is
absolutely continuous.
  \end{itemize}
\end{cor}
\medskip

In particular:
\medskip

\begin{cor} $($Adouani-Marzougui's theorem $($\cite{am12}, Theorem B$)$$)$\label{c:16}
Let $f$  satisfies the assumptions of the main theorem. Then:
\begin{itemize}
    \item [(i)] If $f$ does not have the $(D)$-property then the invariant measure $\mu_{f}$ is singular with respect
to the Lebesgue measure $m$.
    \item [(ii)] If $f$ has the $(D)$-property and $D^{2}f \in
L^{p}(S^{1})$ for some $p>1$, then the measure $\mu_{f}$ is
equivalent to the Lebesgue measure $m$.
  \end{itemize}
\end{cor}
\medskip

\begin{cor} $($\cite{am12}, Corollary 1.5$)$ \label{c:17} Let $f\in \mathrm{PL}(S^{1})$
have irrational rotation number $\alpha$ of bounded type. Then the following are equivalent:
\begin{itemize}
 \item [(i)] $f$ has the ($D$)-property

 \item [(ii)] The measure $\mu_{f}$ is equivalent to the Lebesgue measure $m$.
\end{itemize}
\end{cor}
\medskip

\textbf{Remark.} The main Theorem and in particular Corollary \ref{c:17}  cannot be extended to rotation number not of bounded type, since very recently,
Teplinsky \cite{aT15} constructs an example of a (PL) circle homeomorphism with $4$ non-trivial break points
lying on different orbits that has invariant measure equivalent to the Lebesgue measure $m$.
The rotation number for such example can be chosen either Diophantine or Liouvillean, but not of bounded type.

\section{\bf Notations and preliminary results}
\bigskip

 \textit{\bf 2.1. Dynamical partitions}. Let  $f$ be a homeomorphism of  $S^{1}$  with irrational rotation number  $\alpha = \rho(f)$.
We identify $\alpha$ to its lift $\widehat{\alpha}$ in $]0, 1[$. Let
$(a_{n})_{n\in \mathbb{N}^{\ast}}$ be the partial quotients of
$\alpha$ in the continued fractions expansion. For $n\in
\mathbb{N}^{\ast}$, the fractions $[a_{1},a_{2},\dots,a_{n}]$ are
written in the form of irreducible fractions $\frac{p_{n}}{q_{n}}$.
The sequence $\frac{p_{n}}{q_{n}}$ converges to $\alpha$ and we say
that $\frac{p_{n}}{q_{n}}$ are \textit{rational approximations} of
$\alpha$. Their denominators $q_{n}$ satisfy the following recursion
relation:
$$q_{n} = a_{n}q_{n-1}+q_{n-2}, ~n\geq 2, ~q_{0}=1, ~q_{1}=a_{1}.$$

Let $x_{0}\in S^{1}$ fixed. Denote by:

$$\Delta_{0}^{(n)}(x_{0}) = \begin{cases}
{[x_{0}, f^{q_{n}}(x_{0})]}, &\text{ if  $n$  is  even }  \\
{[f^{q_{n}}(x_{0}), x_{0}]}, &\text{ if  $n$  is  odd }
\end{cases}$$
\
\\
$$\Delta_{i}^{(n)}(x_{0}) : \ = f^{i}\left(\Delta_{0}^{(n)}(x_{0})\right), \ i\in \mathbb{Z}$$
\medskip

\textit{In all the sequel, we deal with the case $n$ odd} (the case
$n$ even is obtained by reversing the orientation of $S^{1}$).
\smallskip

We have then:
\medskip

\begin{lem}[See \cite{yS94}]\label{l:23} The segments  $\Delta_{i}^{(n-1)}(x_{0}) = f^{i}\left(\Delta_{0}^{(n-1)}(x_{0})\right), \ 0\leq
i<q_{n}$  and  $\Delta_{j}^{(n)}(x_{0}) =
f^{j}\left(\Delta_{0}^{(n)}(x_{0})\right), \ 0\leq j<q_{n-1}$ cover
$S^{1}$ and that their interiors are mutually disjoint.
\end{lem}
\medskip

The partition denoted by $\xi_{n}(x_{0}): =
\big\{\Delta_{i}^{(n-1)}(x_{0}); \ \Delta_{j}^{(n)}(x_{0}), \ 0\leq
i<q_{n}, \ 0\leq j<q_{n-1}\big\}$  is called the \textit{$n$-th
dynamical partition} of the point $x_{0}$. It is defined by the
order of points $x_{0}, \dots, f^{q_{n-1}+q_{n}-1}(x_{0})$. The
process of passing from $\xi_{n}(x_{0})$ to $\xi_{n+1}(x_{0})$ is
described by remaining intact the elements
 $\Delta_{j}^{(n)}(x_{0})$, while each element $\Delta_{i}^{(n-1)}(x_{0})$, ($0\leq i<q_{n}$) is partitioned into $a_{n+1}+1$ sub-segments; since  $q_{n+1} = a_{n+1}q_{n}+q_{n-1}$:
$$\Delta_{i}^{(n-1)}(x_{0}) = \Delta_{i}^{(n+1)}(x_{0})~\cup~\bigcup_{s=0}^{a_{n+1}-1}
 \Delta_{i+q_{n-1}+sq_{n}}^{(n)}(x_{0})$$
For every  $n\geq 1,~\xi_{n+1}(x_{0})$  is finer than
$\xi_{n}(x_{0})$ in the following sense: every element
$\xi_{n}(x_{0})$ is a union of elements  $\xi_{n+1}(x_{0})$.
\medskip

The following lemmas are easy to check.
\begin{lem}\label{l:24} Let  $n\geq 1$  be an integer and $y_{1}\in
S^{1}$. Set $y_{2} = f^{q_{n-1}}(y_{1})$,  $y_{3} =
f^{q_{n-1}}(y_{2})$, $\Delta_{i}^{(n-1)}(y_{1}) = f^{i}([y_{1},
y_{2}])$  and  $\Delta_{j}^{(n-1)}(y_{2}) = f^{j}([y_{2}, y_{3}])$,
$0\leq i, \ j < q_{n}$. Then $\Delta_{i}^{(n-1)}(y_{1})\cap
\Delta_{j}^{(n-1)}(y_{2})$ has non empty interior if and only if
\textit{one} of the following conditions holds:
\medskip

\item \begin{enumerate}
\item  $j = i-q_{n-1}$ and $q_{n-1}\leq i<q_{n}$; in this case we have
$$\Delta_{i}^{(n-1)}(y_{1}) = \Delta_{j}^{(n-1)}(y_{2}).$$
\item  $j = q_{n}-q_{n-1}+i$ and $0\leq i<q_{n-1}$; in this case we have

$$\Delta_{j}^{(n-1)}(y_{2}) = f^{q_{n}}(\Delta_{i}^{(n-1)}(y_{1})).$$
\end{enumerate}
\end{lem}
\medskip
\bigskip


\begin{lem}\label{l:25} Let  $n\geq 1$  be an integer and  $y_{1}\in
S^{1}$. Set  $y_{2} = f^{q_{n-1}}(y_{1})$,  $y_{3} =
f^{q_{n-1}}(y_{2})$. Then for every  $x\in S^{1}\backslash
O_{f}(y_{1})$, there exists a unique $0\leq i = i_{n}(x) <q_{n}$
such that $x\in f^{i}(]y_{1},y_{2}[)~\cup
f^{\varphi_{n}(i)}(]y_{2},y_{3}[)$ where $\varphi_{n}:
[0,q_{n}[\longrightarrow [0,q_{n}[$  is a bijective map given by:

$$\varphi_{n}(i) =
\begin{cases}
q_{n}-q_{n-1}+i,~ &\hbox{if} \ \ ~~0\leq i<q_{n-1}\\
i-q_{n-1},~ &\hbox{if} \ \ ~~q_{n-1}\leq i<q_{n}
\end{cases}$$
\end{lem}
\medskip
\bigskip

\textit{\bf 2.2. Useful inequalities.} The proof of the following
results are based on the dynamical partition of  $S^{1}$.
\bigskip

{\it Finzi distortion Lemma}~\cite{aF50}. Let  $f\in
\mathcal{P}(S^{1})$  with irrational rotation number. For all  $z\in
S^{1}$;  $x,~y\in [z,~f^{q_{n-1}}(z)]$ and for all integer $0\leq k
\leq q_{n}$, we have:
$$~e^{-V}\leq\frac{\textrm{Df}^{k}(x)}{\textrm{Df}^{k}(y)}\leq
e^{V}.$$

\
\\
Notice that Finzi distortion Lemma is used for the proof of
Herman-Katok's theorem.
\medskip

The following Lemma plays a key role for studying metrical
properties of the homeomorphism $f$:
\bigskip

{\it Denjoy's inequality}. Let $f\in \mathcal{P}(S^{1})$ with
irrational rotation number. For all $x\in S^{1}$  and all $n\in
\mathbb{N}$, we have  $$e^{-V}\leq \textrm{Df}^{q_{n}}(x)\leq
e^{V}$$ where $V= \textrm{Var}(\log \textrm{Df})$.
\medskip

The proof of this inequality results from Finzi distortion Lemma. It
can be also obtained by applying Denjoy-Koksma's inequality to the
fonction  $\log \textrm{Df}$ (cf. Herman \cite{mH79}).
\medskip

Here and further on, by $m([a,b])$ we mean the smallest length of the two
arcs $[a,b]$ and $[b,a]$ on $S^{1}$.

Using Denjoy's inequality, one has:
\medskip

\begin{cor}[Geometric inequality]\label{c:24} There exists a constant $C>0$ such that for any $x_{0}\in S^{1}$, $n\geq 1$ and any element
$\Delta^{(n)}$ of the dynamical partition $\xi_{n}(x_{0})$, we have
\; $m(\Delta^{(n)})\leq C\lambda^{n}$, where $\lambda =
(1+e^{-V})^{-\frac{1}{2}}< 1$.
\end{cor}
\medskip

From Corollary \ref{c:24} it follows that every orbit of every $x\in
S^{1}$ is dense in $S^{1}$ and this implies the following
generalization of the classical Denjoy theorem:
\medskip

\begin{cor}[Denjoy's theorem: the class P] Let $f\in \mathcal{P}(S^{1})$ with irrational rotation number $\alpha = \rho(f)$.
Then $f$ is topologically conjugate to the  rotation $R_{\alpha}$.
\end{cor}
\medskip

In the following Lemma we have to compare the lengths of iterates
different of intervals.
\medskip

\begin{lem}\label{l:31} Let $f\in \mathcal{P}(S^{1})$ with irrational rotation number $\alpha = \rho(f)$.
Let  $n\in \mathbb{N}^{*}$  and  $z_{1}\in S^{1}$. Set $z_{2} =
f^{q_{n-1}}(z_{1}), \ z_{3} = f^{q_{n-1}}(z_{2})$. Then for any
segments  $K_{1}, \ K_{2}\subset [z_{1},z_{3}]$, one has:\\

$$e^{-2V}\dfrac{m(K_{1})}{m(K_{2})}\leq
\dfrac{m(f^{j}(K_{1}))}{m(f^{j}(K_{2}))}\leq
e^{2V}\dfrac{m(K_{1})}{m(K_{2})}$$ for all \; $j=-q_{n}, \dots
,0,\dots, q_{n}$.
\end{lem}
\medskip

\begin{proof} If  $j= q_{n}$,  Lemma \ref{l:31} is a consequence of Denjoy's inequality.
We suppose that  $0\leq j<q_{n}$. Let  $\mathcal{L}$ be the set of
segments in  $[z_{1},z_{2}]$
 or in  $[z_{2},z_{3}]$. If  $K\in \mathcal{L}$ then  $K\subset
[z,f^{q_{n-1}}(z)]$  with  $z=z_{1}~\textrm{or}~z=z_{2}$.

 Since $$\frac{m(f^{j}(K))}{\textrm{Df}^{j}(z_{2})} = \int_{K}\frac{\textrm{Df}^{j}(t)}{\textrm{Df}^{j}(z_{2})}dm(t),$$
\medskip

by Finzi distortion Lemma applied to $f$  on $[z,f^{q_{n-1}}(z)]$,
we get:

   $$e^{-V}\leq \frac{\textrm{Df}^{j}(t)}{\textrm{Df}^{j}(z_{2})}\leq
e^{V}, \ \textrm{ for every } \ t\in K.$$
\medskip

By integrate this last inequality, we obtain
\begin{equation}e^{-V}m(K)\leq
\frac{m(f^{j}(K))}{\textrm{Df}^{j}(z_{2})}\leq e^{V}m(K)
\end{equation}
\medskip

Let  $K_{1}, \ K_{2}\subset [z_{1},z_{3}]$  be two segments. We
distinguish three cases:
\medskip

a) $K_{1}\in \mathcal{L}$  and  $K_{2}\in \mathcal{L}$: in this
case, the inequality (2.1) implies: \\
$$e^{-V}m(K_{1})\leq
\frac{m(f^{j}(K_{1}))}{\textrm{Df}^{j}(z_{2})}\leq
e^{V}m(K_{1})$$~~and~~
$$e^{-V}m(K_{2})\leq \frac{m(f^{j}(K_{2}))}{\textrm{Df}^{j}(z_{2})}\leq
e^{V}m(K_{2}),$$ \\
therefore
$$e^{-2V}\frac{m(K_{1})}{m(K_{2})}\leq\frac{m(f^{j}(K_{1}))}{m(f^{j}(K_{2}))}
 \leq e^{2V}\frac{m(K_{1})}{m(K_{2})}.$$
\medskip

b) $K_{1}\notin \mathcal{L}$ and $K_{2}\in \mathcal{L}$: in this
case, take $K_{1} = K_{1}^{\prime}\cup K_{1}^{\prime\prime}$  where
$K_{1}^{\prime}\subset [z_{1},z_{2}]$ and
$K_{1}^{\prime\prime}\subset [z_{2},z_{4}]$ are two segments, then
$K_{1}^{\prime}\in \mathcal{L}$ and
$K_{1}^{\prime\prime}\in\mathcal{L}$. By ~ a), we see that
\\
$$e^{-2V}\frac{m(K_{1}^{\prime})}{m(K_{2})}\leq\frac{m(f^{j}(K_{1}^{\prime}))}{m(f^{j}(K_{2}))}
 \leq e^{2V}\frac{m(K_{1}^{\prime})}{m(K_{2})}$$ and
 $$e^{-2V}\frac{m(K_{1}^{\prime\prime})}{m(K_{2})}\leq\frac{m(f^{j}(K_{1}^{\prime\prime}))}{m(f^{j}(K_{2}))}
 \leq e^{2V}\frac{m(K_{1}^{\prime\prime})}{m(K_{2})}.$$
\medskip

Adding these inequalities, we get
  $$e^{-2V}\frac{m(K_{1})}{m(K_{2})}\leq\frac{m(f^{j}(K_{1}))}{m(f^{j}(K_{2}))}
 \leq e^{2V}\frac{m(K_{1})}{m(K_{2})}.$$
\medskip

c) $K_{1}\notin \mathcal{L}$  and  $K_{2}\notin \mathcal{L}$: in
this case, take  $K_{2} = K_{2}^{\prime}\cup K_{2}^{\prime\prime}$
where  $K_{2}^{\prime}\subset [z_{1},z_{2}]$ and
$K_{2}^{\prime\prime}\subset [z_{2},z_{3}]$ are two segments. By
applying ~ b) to  $K_{1}$  and
 $K_{2}^{\prime}$ (resp. $K_{2}^{\prime\prime}$), we obtain the inequality of the Lemma.
Apply the same result to $f^{-1}$ instead of $f$ (since the
convergent sequences of the continued fractions of $\varrho(f^{-1})=1-\varrho(f)$
and $\varrho(f)$ have the same denominators).
\end{proof}
\bigskip

\textit{\bf 2.3. Some ratio tools.} Let us introduce the notion of
ratio distortion with respect to a homeomorphism.
\medskip

\begin{defn} Let  $x_{1},x_{2},x_{3}$  be real numbers such that  $x_{1}<x_{2}<x_{3}$.
\medskip
\
\\
i)  The \textit{ratio} of the triple ($x_{1},x_{2},x_{3}$) is the
real number defined as

 $$\textrm{r}(x_{1},x_{2},x_{3}) =
\frac{x_{2}-x_{1}}{x_{3}-x_{2}}.$$

\
\\
ii) The \textit{ratio distortion} of the triple
($x_{1},x_{2},x_{3}$) with respect to a strictly increasing function
$F: \mathbb{R}\longrightarrow \mathbb{R}$  is the real number
defined as

$$\textrm{Dr}_{F}(x_{1},x_{2},x_{3}) =
\frac{\textrm{r}(F(x_{1}),F(x_{2}),F(x_{3}))}{r(x_{1},x_{2},x_{3})}
= \frac{F(x_{2})-F(x_{1})}{x_{2}-x_{1}}
\frac{x_{3}-x_{2}}{F(x_{3})-F(x_{2})}$$
\end{defn}
\medskip

Notice that:
\medskip

- The ratio $\textrm{r}$ is translation invariant:
$$\textrm{r}(x_{1}+a,x_{2}+a,x_{3}+a) = \textrm{r}(x_{1},x_{2},x_{3}), \ a\in \mathbb{R}.$$

- The ratio distortion $\textrm{Dr}$ is multiplicative with respect
to composition: for two functions $F$ and $G$ on $\mathbb{R}$, we
have
$$\textrm{Dr}_{G\circ F}(x_{1},x_{2},x_{3}) = \textrm{Dr}_{G}(F(x_{1}),F(x_{2}),F(x_{3}))\textrm{Dr}_{F}(x_{1},x_{2},x_{3}).$$
\medskip

Let $z_{1},z_{2},z_{3}\in S^{1}$ with lifts $\widehat{z_{1}}, \
\widehat{z_{2}}, \ \widehat{z_{3}}$ and $z_{1}\prec z_{2}\prec
z_{3}\prec z_{1}$  ordered on $S^{1}$. Define

$$\widetilde{z_{i}}:=
\begin{cases}
    \widehat{z_{i}}, & \ \textrm{ if } \  \widehat{z_{1}}< \widehat{z_{i}} <1\\
    1+\widehat{z_{i}}, & \ \textrm{ if } \  0\leq \widehat{z_{i}}< \widehat{z_{1}}
 \end{cases}$$
\
\\
where $i=2,3$. Obviously we have $\widetilde{z_{1}} <
\widetilde{z_{2}} < \widetilde{z_{3}}$. The vector
$(\widetilde{z_{1}},\widetilde{z_{2}},\widetilde{z_{3}}) \in
\mathbb{R}^{3}$ is called the\textit{ lifted vector } of
$(z_{1},z_{2}, z_{3})\in (S^{1})^{3}$.

Consider now a homeomorphism $f$ of $S^{1}$  with lift
$\widehat{f}$. We define the \textit{ratio distortion} of the triple
$(z_{1},z_{2},z_{3})$, $z_{1},z_{2},z_{3}\in S^{1}$ with respect to
$f$  by
$$\textrm{Dr}_{f}(z_{1},z_{2},z_{3}): = \textrm{Dr}_{\widehat{f}}(\widetilde{z_{1}},\widetilde{z_{2}},\widetilde{z_{3}})$$  where
$(\widetilde{z_{1}},\widetilde{z_{2}},\widetilde{z_{3}})$ is the
lifted vector of $(z_{1},z_{2}, z_{3})$. This definition is
independent of the lift of $f$  since the ratio is translation invariant.
\medskip

\begin{defn}\label{d:32} Let $R>1$ be a real number and $x_{0}\in S^{1}$. We say that a triple
$(z_{1},z_{2},z_{3})$~ of~ $S^{1}$ ($z_{1}\prec z_{2}\prec z_{3})$
 satisfies the conditions $(a)$~~\textrm{and}~~$(b)$ for the point $x_{0}$ and
the constant $R$ if:\\
\medskip

\begin{itemize}
 \item [(a)]: $R^{-1}\leq\dfrac{m([z_{2},~z_{3}])}{m([z_{1},~z_{2}])}\leq R$
\
\\
\item [(b)]: $\underset{1\leq i\leq 3}{\max}~m([x_{0},z_{i}])\leq R m([z_{1},z_{2}])$
\end{itemize}
\end{defn}
\medskip

We call two intervals in $S^{1}$ $R$-\textit{comparable} if the
ratio of their lenghts is in $[R^{-1}, R]$.
\bigskip
\medskip

\textit{\bf 2.4. Reduction.} In this subsection, we will reduce any
homeomorphism $f\in \mathcal{P}(S^{1})$ with several break (or
non-break)
 points to the one with break points lying on different orbits.
\medskip

\begin{defn} Let $c\in C(f)$. A \emph{maximal $f$-connection} of  $c$ is a segment
\ $[f^{-p}(c),\dots, f^{q}(c)]:= \{f^{s}(c): \ -p \leq s\leq q\}$
of the orbit  $O_{f}(c)$  which contains all the break points of $f$
contained on  $O_{f}(c)$  and such that  $f^{-p}(c)$  (resp.
$f^{q}(c)$)  is the first (resp. last) break point of  $f$  on
$O_{f}(c)$.
\end{defn}
\medskip
\
\\
 We have the following properties:
\medskip
\
\\
- Two break points of $f$ are on the same maximal $f$-connection, if
and only if, they are on the same orbit.

- Two distinct maximal $f$-connections are disjoint.
\bigskip

Denote by \
\\
\begin{itemize}

 \item [-]  $ M_{i}(f) = [c_{i},\dots, f^{N_{i}}(c_{i})], \ (N_{i}\in \mathbb{N}),$ the maximal $f$-connections of
$c_{i}\in C(f)$, ($0\leq i\leq p$).

 \item [-]  M$(f) = \coprod _{i=0}^{p}M_{i}(f)$.
\end{itemize}
\medskip
So, we have the decomposition: $C(f) = \coprod_{i=0}^{p} C_{i}(f)$
where,  $C_{i}(f) = C(f)\cap M_{i}(f), \ 0\leq i\leq p$. We also
have $$\underset{d\in C_{i}(f)}\prod \sigma_{f}(d) = \underset{d\in
M_{i}(f)}\prod \sigma_{f}(d).$$
\medskip
\
\\
\begin{prop}[\cite{am15}, Theorem 2.1]\label{p:623}
Let $f\in \mathcal{P}(S^{1})$ with irrational rotation number, and
let $\big(k_{0},\dots,k_{p}\big)\in \mathbb{Z}^{p+1}$. Then there
exists a piecewise quadratic homeomorphism $K\in \mathcal{P}(S^{1})$
such that $F:= K \circ f \circ K^{-1}\in \mathcal{P}(S^{1})$ with
$C(F) \subset \{K(f^{k_{i}}(c_{i}))= F^{k_{i}}(K(c_{i})); i=0,1,
\dots, p\}$ and such that $\sigma_{F}(F^{k_{i}}(K(c_{i}))) =
\pi_{s,O_{f}(c_{i})}(f),~i=0,1, \dots, p$.
\end{prop}
\medskip

\begin{cor}[\cite{am15}, Corollary 2.3]\label{c:21} Let $f\in \mathcal{P}(S^{1})$ with irrational rotation number.
Suppose that $f$ satisfies the (KO) condition. Then, there exists a
piecewise quadratic homeomorphism $K\in \mathcal{P}(S^{1})$ such
that: $F = K\circ f \circ K^{-1}\in \mathcal{P}(S^{1})$ with $C(F)\subset\{K(c_{0}),\dots,K(c_{p})\}$, where $c_{0},\dots, c_{p}\in C(f)$
are on pairwise distinct orbits.
\end{cor}
\medskip


In particular, we have:
\medskip

\begin{cor}\label{c:22} Let $f$, $K$ and $F$ as in Corollary \ref{c:21}.
\begin{itemize}
\item [(i)] If $f$ does not satisfy the ($D$)-property, then there exists
$0\leq i\leq p$  such that $K(c_{i})$ is the unique break point of
$F$ in its orbit.
 \item [(ii)] If $f$ satisfies the ($D$)-property then $F$ is a $C^{1}$-diffeomorphism with $\mathrm{DF}$ absolutely continuous on
$S^{1}$.
\end{itemize}
\end{cor}
\medskip
\bigskip

\begin{prop}\label{p:23} Let $f, ~g\in \mathcal{P}(S^{1})$ and have the same irrational rotation number such that the break points of
$f$ (resp. $g$) belong to pairwise distinct $f$-orbits (resp.
$g$-orbits). Let $h$ be the conjugating homeomorphism from $f$ to
$g$. Then there exist an integer $q\geq p$,
$c_{0},c_{1},\dots,c_{q}; \ d_{0},d_{1},\dots,d_{q}$ belong to pairwise distinct orbits $f$-orbits and $g$-orbits respectively
such that $C(f)\subset \{c_{0},c_{1},\dots,c_{q}\}$,
$C(g)\subset \{d_{0},d_{1},\dots,d_{q}\}$, and a homeomorphism $u$
of $S^{1}$ such that $G: = u \circ f \circ u^{-1}\in
\mathcal{P}(S^{1})$ has break points belong to pairwise distinct
$G$-orbits and satisfies
\\
\begin{itemize}

 \item [(i)]  $C(G)\subset \{u(c_{i}): ~i=0,1, \dots, q\}$.
 \item [(ii)] $\sigma_{G}(u(c_{i})) = \sigma_{g}(d_{i})~,~i=0,1, \dots, q$.
 \item [(iii)] $\pi_{s}(G)=\pi_{s}(g)$.
 \item [(iii)] $h$ is singular if and only if so is $u$.
\end{itemize}
\end{prop}
\medskip
\
\\
\textit{From now on we denote by \ \ $B:=  \{c_{0},\dots, c_{q}\}$} and $c_{0}=c$.
\medskip

In the sequel, we may assume that $f, ~g\in \mathcal{P}(S^{1})$ have
the same irrational rotation number $\alpha$ and satisfy the
following:
\medskip

- The points of $B$ belong to pairwise distinct $f-$ orbits.

- The break points of $g$ belong to pairwise distinct $g$-orbits.

- $C(f)\subset B$.

- The maps $f, ~g$ satisfy the Katznelson-Ornstein (KO) condition.

- The conjugating $h$ is such that $C(g)\subset h(B)$.
\bigskip
\medskip

\section{\bf Primary Cells}
\medskip
\medskip

Let  $x_{0}\in S^{1}$ and $f\in \mathcal{P}(S^{1})$
with irrational rotation number $\alpha =\rho(f)$. Let $c\in B$
and $n$ an odd integer. By Lemma \ref{l:23}, either $c\in
\Delta^{(n-1)}_{i_{n}(c)}(x_{0})$ for some $0\leq i_{n}(c) < q_{n}$
or $c\in\Delta^{(n)}_{i_{n}(c)}(x_{0})$ for some $0\leq i_{n}(c)<
q_{n-1}$. Set

\[\ y_{2} = f^{-i_{n}(c)}(c), \ y_{1} =
f^{-q_{n-1}}(y_{2}), \ y_{3} = f^{q_{n-1}}(y_{2}).\]
\smallskip

Notice that $y_{1},y_{2}$ and $y_{3}$ are defined with respect to
$f,x_{0},c$ and the number $i_{n}(c)$ depends on $x_{0}$.
\medskip

Let $\delta > 0$ and $U_{\delta}(x_{0})$ a
$\delta$-neighbourhood of $x_{0}$.
\medskip

\begin{prop}$($cf. \cite{am12}, Proposition 3.1$)$ \label{p:41} Under the notations above, there exists $N = N(x_{0}, \delta)\in \mathbb{N}$ such that for all $n\geq N$, there is
a triple $(y_{1},y_{2},y_{3})=(y_{1}(n),y_{2}(n),y_{3}(n))_{n\geq
N}$ with the following properties:
\medskip

 \begin{itemize}
 \item [(c-0)] $\left(y_{1},y_{2},y_{3}\right) \left(\mathrm{resp}.
 ~ (f^{q_{n}}\left(y_{1}),f^{q_{n}}(y_{2}),f^{q_{n}}(y_{3}\right)\right))_{n\geq N}\subset
U_{\delta}(x_{0})$\\

\item [(c-1)] $~y_{2} \in \Delta^{(n-1)}_{0}(x_{0})~\mathrm{ or }~y_{2}\in \Delta^{(n)}_{0}(x_{0})$\\

 \item [(c-2)] $~m\left(f^{j}([y_{1},y_{3}])\right)\leq K\lambda^{n}$, for every $0\leq j<q_{n}$, where $\lambda = (1+e^{-V})^{-\frac{1}{2}}$  and  $K$ a constant independent of  $n$.

 \item [(c-3)] $~\sum_{j=0}^{q_{n}-1}m\left(f^{j}([y_{1},y_{3}])\right)\leq 2$\\

 \item [(c-4)] The triples  $\left(y_{1},y_{2},y_{3}\right)$~ and
$(f^{q_{n}}(y_{1}),f^{q_{n}}(y_{2}),f^{q_{n}}(y_{3}))$ satisfy
conditions $(a)$ and $(b)$ of definition \ref{d:32} for
$x_{0}$ and the constant $R = e^{3V}+ e^{V}+1$.\\

 \item [(c-5)] \begin{itemize}
\item [-] $f^{i}([y_{1},~y_{3}])$ does not contain $c$ for every $0\leq i\neq i_{n}(c)<q_{n}$.
\item [-] $f^{i_{n}(c)}([y_{1},y_{3}])$ does not contain any break point of $f$ other than $c$.
 \end{itemize}

\end{itemize}
\end{prop}
\medskip

\textbf{Definition $($primary cell$)$}. The triple
$\left(y_{1},y_{2},y_{3}\right)$ given in Proposition \ref{p:41} is
called a \textit{primary cell associated to  $(f,x_{0},c,\delta,N)$}
for the constant $R = e^{3V}+ e^{V}+1$.
\bigskip

In the sequel, set $\mathbb{N}_{0}=\{n \in \mathbb{N}~,~n \geq N\}$.

Let $c \in B$, $d\notin O_{f}(c)$ and
$\left(y_{1},y_{2},y_{3}\right): =
(y_{1}(n),y_{2}(n),y_{3}(n))_{n\in \mathbb{N}_{0}}$ be the primary cell
associated to  $(f,x_{0},c,\delta,N)$. By Lemma \ref{l:25} and for
every and $n\in \mathbb{N}_{0}$, denote by

$$t_{n}(d)   = \dfrac{m([f^{-i_{n}(d)}(d), y_{2}])}{m([f^{j_{n}(d)-i_{n}(d)}(y_{2}), y_{2}])}$$
\\
where \;$j_{n}(d) =
\varphi_{n}(i_{n}(d))$ and $i_{n}(d)$ is the unique integer $0 \leq
i_{n}(d)<q_{n}$ such that $d \in f^{i_{n}(d)}([y_{1}, y_{2}])$ or $d
\in f^{j_{n}(d)}([y_{1},y_{2}])$.
\bigskip

One always has $f^{-i_{n}(d)}(d)\in [f^{q_{n}}y_{1}, y_{2}]$.
Indeed, by Lemma \ref{l:25}, there exists a unique $0\leq i_{n}(d) <
q_{n}$ such that either
 $d\in f^{i_{n}(d)}([y_{1}, y_{2}])$ or $d\in f^{j_{n}(d)}([y_{2}, y_{3}])$. In the first case, $f^{-i_{n}(d)}(d)\in [y_{1}, y_{2}]$ and in the second case, $f^{-i_{n}(d)}(d)\in [f^{q_{n}}y_{1}, f^{q_{n}}y_{2}]\subset [f^{q_{n}}y_{1}, y_{2}]$.
One can write $$t_{n}(d)   = \left\{
\begin{array}{ll}
       \dfrac{m([f^{-i_{n}(d)}(d), y_{2}])}{m([y_{1}, y_{2}])}, & \hbox{if }q_{n-1}\leq i_{n}(d)< q_{n} \\\\
    \dfrac{m([f^{-i_{n}(d)}(d), y_{2}])}{m([f^{q_{n}}y_{1}, y_{2}])}, & \hbox{if } 0\leq i_{n}(d)< q_{n-1}
\end{array}
\right.
$$
\medskip

Let $d\in B\backslash\{c\}$ and $A$ be an infinite subset of $\mathbb{N}_{0}$. Denote by
\medskip

\textbullet \; $T(d, A)$ the set of limit points of $(t_{n}(d))_{n\in A}$. Therefore
 $T(d, A)\subset [0,1]$.
\medskip

We say that $d$ satisfies:
\medskip

\textbullet \; The $P_{1}(d,A)$-property if $T(d, A)\cap ~]0,1[\neq \emptyset$
and for every $n\in A, ~0 \leq i_{n}(d) <q_{n-1}$.

\textbullet \; The $P_{2}(d,A)$-property if $T(d, A)\cap ~]0,1[\neq \emptyset$
and for every $n \in A, ~q_{n-1} \leq i_{n}(d) <q_{n}$.
\bigskip
\medskip

Set \begin{align*} Q_{11}( A) & = \{d\in B\backslash \{c\}:
~\textrm{d ~satisfies ~the } P_{1}(d, A)-\textrm{property}
\}\\
Q_{12}(A) & = \{d\in B\backslash \{c\}:  ~\textrm{d ~satisfies ~the
} P_{2}(d, A)-\textrm{property}
\}\\
Q_{1}(A) & = Q_{11}(A)\cup Q_{12}(A)\\
Q_{2}(A) & = \{d\in B\backslash \{c\}:  T(d, A) = \{0\}\}\\
Q_{3}( A) & = \{d\in B\backslash \{c\}: T(d, A) = \{1\}\}\\
Q_{4}( A) & = \{d\in B\backslash \{c\}: T(d, A) = \{0,1\}\}
\end{align*}
\bigskip

\begin{prop}\label{p:461}
Assume that the rotation number $\rho(f)$ of $f$ is of bounded type. Then for every points $c, c_{1}$ in
$B(f)$, $c_{1}\notin O_{f}(c)$, there exists an infinite subset $\mathbb{M}_{c,c_{1}}$ of
$\mathbb{N}_{0}$ such that $c_{1}\in Q_{1}(\mathbb{M}_{c,c_{1}})$.
\end{prop}
\medskip

We need two lemmas.
\medskip

\begin{lem}\label{l:4611}
Let $0\leq k<q_{n+1}+q_{n}$ be an integer. Then
$\Delta_{0}^{(n-1)}(c)\cap \Delta^{(n)}_{k}(c)$ has non empty
interior if and only if $k=q_{n-1}+q_{n+1}$ or $k = q_{n-1}+jq_{n}$ for some $0\leq j\leq
a_{n+1}$.
\end{lem}
\smallskip

\begin{proof}
- Assume that $k=q_{n-1}+q_{n+1}$ or $k = q_{n-1}+jq_{n}$ for some $0\leq j \leq
a_{n+1}$. Since $$\Delta_{0}^{(n-1)}(c)= \Delta_{0}^{(n+1)}(c) \coprod_{0 \leq
j<a_{n+1} }~~\Delta^{(n)}_{q_{n-1}+jq_{n}}(c),$$ so we have
$$\Delta_{0}^{(n-1)}(c)\cap \Delta^{(n)}_{k}(c) = \begin{cases}
\Delta^{(n)}_{k}(c), & \textrm{if } k=q_{n-1}+jq_{n} \ \textrm{ with } 0\leq j <a_{n+1}, \\
 [c,f^{q_{n+1}}(c)], & \textrm{if } k=q_{n+1}, \\
 [f^{q_{n-1}+q_{n}+q_{n+1}}(c),f^{q_{n-1}}(c)], & \textrm{if } k=q_{n-1}+q_{n+1}.
\end{cases}$$

In either case, $\Delta_{0}^{(n-1)}(c)\cap \Delta^{(n)}_{k}(c)$ has non empty
interior. Conversely we distinguish two cases:

- $0 \leq k< q_{n+1}$: if $k \notin
\{q_{n-1}+jq_{n}:~0\leq j< a_{n+1}\}$, then by Lemma
\ref{l:23}, $\Delta_{0}^{(n-1)}(c)\cap \Delta^{(n)}_{k}(c)$ has
empty interior.

- $q_{n+1}\leq k< q_{n+1}+q_{n}$: If $k\notin \{q_{n+1}, q_{n-1}+q_{n+1}\}$, then $k=l+q_{n+1}$ with $1\leq l< q_{n}$ and 
$l\neq q_{n-1}$. Thus we have
\begin{align*}
f^{-q_{n+1}}\left(\Delta_{0}^{(n-1)}(c)\cap
\Delta^{(n)}_{k}(c)\right)=[f^{-q_{n+1}}(c),~f^{-q_{n+1}+q_{n-1}}(c)]\cap
\Delta^{(n)}_{l}(c)\\
 \subset \left(\Delta_{0}^{(n)}(c)\cup
\Delta_{0}^{(n-1)}(c)\right)\cap \Delta^{(n)}_{l}(c)
\end{align*}
Hence $\Delta_{0}^{(n-1)}(c)\cap \Delta^{(n)}_{k}(c)$ has empty
interior.
\end{proof}
\medskip

\begin{lem}\label{l:4612}
Let $0 \leq k < (a_{n+1}+2)q_{n}$ be an integer. Then
$\Delta^{(n-1)}_{k}(c)\cap \Delta_{0}^{(n)}(c)$ has non empty
interior if and only if $k = q_{n}+q_{n+1}$ or $k = jq_{n}$ for some
$1\leq j\leq a_{n+1}+1$.
\end{lem}
\medskip

\begin{proof} We have $$\Delta_{0}^{(n)}(c)\cap \Delta^{(n-1)}_{k}(c) = \begin{cases}
\Delta_{0}^{(n)}(c), & \textrm{if }  k=jq_{n} \textrm{ with } \  1 \leq j \leq a_{n+1}, \\
[f^{q_{n}}(c),f^{q_{n}+q_{n+1}}(c)], & \textrm{if } k=(a_{n+1}+1)q_{n}, \\
[f^{q_{n}+q_{n+1}}(c),c], & \textrm{if } k= q_{n}+q_{n+1}.
\end{cases}$$
Hence $\Delta_{0}^{(n)}(c)\cap \Delta^{(n-1)}_{k}(c)$ has non empty interior. Conversely, we distinguish three cases:

- $0 \leq k< q_{n}$: In this case, by Lemma \ref{l:23}, $\Delta^{(n-1)}_{k}(c)\cap \Delta_{0}^{(n)}(c)$ has empty
interior.

-  $q_{n}\leq k < (a_{n+1}+1)q_{n}$: If $k\notin \{ jq_{n}: \ 1\leq j \leq a_{n+1}\}$, then $k=l+jq_{n}$ with $1\leq l <q_{n}$ 
and $1\leq j\leq a_{n+1}$. Thus we have 
$$f^{-jq_{n}}\left(\Delta^{(n-1)}_{k}(c)\cap
\Delta_{0}^{(n)}(c)\right)=\Delta^{(n-1)}_{l}(c)\cap
f^{-jq_{n}}(\Delta_{0}^{(n)}(c))\subset \Delta^{(n-1)}_{l}(c)\cap
\Delta_{0}^{(n-1)}(c),$$ so  $\Delta^{(n-1)}_{k}(c)\cap
\Delta_{0}^{(n)}(c)$ has empty interior.

-  $(a_{n+1}+1)q_{n}\leq k < (a_{n+1}+2)q_{n}$: If $k\notin \{(a_{n+1}+1)q_{n}, \ q_{n+1}+q_{n}\}$, then $k=l+(a_{n+1}+1)q_{n}$ 
with $1 \leq l<q_{n}$ and $l\neq
q_{n-1}$. Thus we have
\begin{align*}
f^{-(a_{n+1}+1)q_{n}}\left(\Delta^{(n-1)}_{k}(c)\cap
\Delta_{0}^{(n)}(c)\right) & = \Delta^{(n-1)}_{l}(c)\cap
f^{-(a_{n+1}+1)q_{n}}(\Delta_{0}^{(n)}(c))\\
& \subset \Delta^{(n-1)}_{l}(c)\cap \left(\Delta^{(n-1)}_{0}(c)\cup
f^{q_{n-1}}(\Delta_{0}^{(n-1)}(c))\right).
\end{align*}

 Hence $\Delta^{(n-1)}_{k}(c)\cap \Delta_{0}^{(n)}(c)$ has empty
interior.
\end{proof}
\bigskip

\textit{Proof of Proposition \ref{p:461}}. Suppose that there exist two distinct points $c,c_{1}$ in $B(f)$
such that for any infinite subset $M$ of $\mathbb{N}_{0}$ we have $T(c_{1}, M)\subset \{0,1\}$. In particular for
$M= \mathbb{N}_{0}$. So we have three possible cases.
\medskip

\textbf{Case 1}: $T(c_{1}, \mathbb{N}_{0})=\{0\}$ i.e. $c_{1}\in
Q_{2}(\mathbb{N}_{0})$. Since $$f^{-i_{n}(c_{1})}(c_{1})\in
[f^{q_{n}}(y_{1}),~y_{2}] = \Delta_{0}^{(n)}(y_{1})\cup
\Delta_{0}^{(n-1)}(y_{1}),$$ $$m(\Delta_{0}^{(n-1)}(y_{1}))\leq
m([f^{q_{n}}(y_{1}),~y_{2}])\leq
(1+e^{V})m(\Delta_{0}^{(n-1)}(y_{1}))$$ and $$\lim_{n \to
+\infty}~\frac{m([f^{-i_{n}(c_{1})}(c_{1}),~
y_{2}])}{m([f^{q_{n}}(y_{1}),~y_{2}])}=0,$$ there exists an integer
$n_{0}\in \mathbb{N}$ such that for every $n\geq n_{0},
~f^{-i_{n}(c_{1})}(c_{1})\in \Delta_{0}^{(n-1)}(y_{1})$ and
$$\lim_{n\to +\infty}~\frac{m([f^{-i_{n}(c_{1})}(c_{1}),~y_{2}])}{m(\Delta_{0}^{(n-1)}(y_{1}))}=0.$$
For $n\in \mathbb{N}$, set
\medskip

\[ ~ s_{n}: = i_{n}(c_{1})-i_{n}(c)-q_{n-1}\]

\[l_{n}: = \frac{m(\Delta^{(n)}_{s_{n}}(c))}{m(\Delta^{(n-1)}_{s_{n}}(c))}, ~~~~ \ 
 b_{n}: = \frac{m([c_{1},~f^{s_{n}+q_{n-1}}(c)])}{m(\Delta^{(n-1)}_{s_{n}}(c))}\]
\medskip

Then $-q_{n}<s_{n}+q_{n-1}<q_{n}$ and by (Lemma \ref{l:31} and
Denjoy's inequality), we have for every $n\geq n_{0},~~c_{1}\in
\Delta^{(n-1)}_{s_{n}}(c)$ and
\begin{equation}\label{(3.1)}
\lim_{n\to +\infty} b_{n}=0
\end{equation}
In particular, for every $n\geq n_{0}$, $c_{1}$ is an interior point
of $\Delta^{(n-1)}_{s_{n}}(c)\cap \Delta^{(n)}_{s_{n+1}}(c)$ and
hence $\Delta_{0}^{(n-1)}(c)\cap \Delta^{(n)}_{s_{n+1}-s_{n}}(c)$
has a non empty interior.

\medskip

- If $s_{n+1}-s_{n} \geq 0$, then $s_{n+1}-s_{n} <
q_{n+1}-(q_{n}-q_{n-1})$ and by Lemma \ref{l:4611}
$$s_{n+1}-s_{n}\in \left\{ q_{n-1}+jq_{n}: 0\leq j<
a_{n+1}\right\}$$
\medskip

-If $s_{n+1}-s_{n} \leq 0$, then $s_{n}-s_{n+1} < (a_{n+1}+1)q_{n}$
and by Lemma \ref{l:4612}
$$s_{n+1}-s_{n}\in \{-jq_{n}: \ 1\leq j\leq a_{n+1}\}.$$
\bigskip

Let $n\geq n_{0}$. Since $\Delta^{(n-1)}_{s_{n}}(c)\subset
\bigcup_{0 \leq j \leq
a_{n+1}}~f^{q_{n-1}+jq_{n}}(\Delta^{(n)}_{s_{n}}(c))$, so by
Denjoy's inequality, we have
$$\frac{m(\Delta^{(n-1)}_{s_{n}}(c))}{m(\Delta^{(n)}_{s_{n}}(c))}<
\sum_{0\leq j\leq a_{n+1}}~e^{(1+j)V}$$

So, $$l_{n}^{-1}< \frac{e^{V}}{e^{V}-1}~(e^{(1+a_{n+1})V}-1)$$

Hence,
\begin{equation}\label{(3.2)}
  \textrm{If}~~ \lim_{n\to + \infty }l_{n}=0 \ \textrm{ then } ~ \lim_{n\to + \infty }a_{n+1} = +\infty
\end{equation}
\medskip

We distinguish four sub-cases.
\bigskip

Case 1.1: $s_{n+1}-s_{n} \in \{q_{n-1}+jq_{n}~:~ 1\leq j<
a_{n+1}\}$. Then by the Denjoy's inequality, we have
\begin{align*}
l_{n} & \leq
e^{V}\frac{m(\Delta^{(n)}_{s_{n}+q_{n-1}}(c))}{m(\Delta^{(n-1)}_{s_{n}}(c))}\\
& \leq
e^{V}\frac{m([c_{1},~f^{s_{n}+q_{n-1}}(c)])}{m(\Delta^{(n-1)}_{s_{n}}(c))} \\
& = e^{V}b_{n}
\end{align*}

(since $c_{1} \in \Delta^{(n-1)}_{s_{n}}(c)\backslash
\Delta^{(n)}_{s_{n}+q_{n-1}}(c)$).\\
\medskip

Case 1.2: $s_{n+1}-s_{n}=q_{n-1}$. Then by the Denjoy's inequality, we have
\begin{align*}l_{n}& \leq e^{V}\frac{m(\Delta^{(n)}_{s_{n+1}}(c))}{m(\Delta^{(n-1)}_{s_{n}}(c))} \\
& \leq e^{V}\frac{m([f^{s_{n+1}+q_{n}}(c),~c_{1}])}{m(\Delta^{(n)}_{s_{n+1}}(c))}+
e^{V}\frac{m([c_{1},~f^{s_{n}+q_{n-1}}(c)])}{m(\Delta^{(n-1)}_{s_{n}}(c))} \\
& = e^{V}(b_{n}+b_{n+1})
\end{align*}
(since $c_{1}\in
\Delta^{(n)}_{s_{n+1}}(c) = \Delta^{(n)}_{s_{n}+q_{n-1}}(c)\subset
\Delta^{(n-1)}_{s_{n}}(c)$).\\
\medskip

Case 1.3: $s_{n+1}-s_{n} \in \{-jq_{n}~:~1 \leq j < a_{n+1} \}$. Then by
the Denjoy's inequality, we have
\begin{align*} l_{n}& \leq e^{V}\frac{m(\Delta^{(n)}_{s_{n}+q_{n-1}}(c))}{m(\Delta^{(n-1)}_{s_{n}}(c))}\\
& \leq e^{V}\frac{m([c_{1},~f^{s_{n}+q_{n-1}}(c)])}{m(\Delta^{(n-1)}_{s_{n}}(c))} \\
& = e^{V}b_{n}
\end{align*}
(since $c_{1} \in \Delta^{(n-1)}_{s_{n}}(c)\backslash
\Delta^{(n)}_{s_{n}+q_{n-1}}(c)$ as in Case 1.1).\\
\medskip

Case 1.4: $s_{n+1}-s_{n} = q_{n-1}-q_{n+1}$. Then by the Denjoy's iequality, we have

\begin{align*}
l_{n} & = \frac{m(f^{-q_{n-1}+q_{n+1}}(\Delta^{(n)}_{s_{n+1}}(c)))}{m(\Delta^{(n-1)}_{s_{n}}(c))}\\
& \leq e^{2V}
\frac{m(\Delta^{(n)}_{s_{n+1}}(c))}{m(\Delta^{(n-1)}_{s_{n}}(c))}\\
& \leq e^{2V}\frac{m([f^{s_{n+1}+q_{n}}(c),~c_{1}])}{m(\Delta^{(n-1)}_{s_{n}}(c))}+ 
e^{2V}\frac{m([c_{1},~f^{s_{n+1}}(c)])}{m(\Delta^{(n-1)}_{s_{n}}(c))}\\
& \leq e^{2V}\frac{m([f^{s_{n+1}+q_{n}}(c),~c_{1}])}{m(\Delta^{(n)}_{s_{n+1}}(c))}+
e^{2V}\frac{m([c_{1},~f^{s_{n}+q_{n-1}}(c)])}{m(\Delta^{(n-1)}_{s_{n}}(c))}\\
& \leq e^{2V}(b_{n}+b_{n+1})
\end{align*}
(since $c_{1} \in
\Delta^{(n)}_{s_{n+1}}(c) = f^{-q_{n+1}}(\Delta^{(n)}_{s_{n}+q_{n-1}}(c))
\subset \Delta^{(n-1)}_{s_{n}}(c)$).\\

In either of the four cases above, we conclude that for every $n\geq n_{0}$,~$$l_{n}\leq
e^{2V}(b_{n}+b_{n+1}).$$ So by the assertions (\ref{(3.1)}) and
(\ref{(3.2)}) we obtain that
$\underset{n\to + \infty}\lim ~a_{n+1} = +\infty,$ a contradiction with that $\rho(f)$ is of bounded type.\\
\medskip

\textbf{Case 2}: $ T(c_{1}, \mathbb{N}_{0}) = \{1\}$ i.e. $c_{1}\in
Q_{3}(\mathbb{N}_{0})$. In this case, there
exists $n_{1}\in \mathbb{N}$ such that for all $n\geq
n_{1}, ~f^{-j_{n}(c_{1})}(c_{1})\in \Delta^{(n-1)}(y_{2})$ and
$$\lim_{n\to +\infty}~\frac{m([y_{2},~f^{-j_{n}(c_{1})}(c_{1})])}{m(\Delta^{(n-1)}_{0}(y_{2}))} = 0.$$

For $n\geq n_{1}$, set
\medskip

\[r_{n}: = j_{n}(c_{1})-i_{n}(c)\]

\[l^{\prime}_{n}:= \frac{m(\Delta^{(n)}_{r_{n}}(c))}{m(\Delta^{(n-1)}_{r_{n}}(c))}, ~~~~ \ 
b^{\prime}_{n}:= \frac{m([f^{r_{n}}(c),~c_{1}])}{m(\Delta^{(n-1)}_{r_{n}}(c))}\]
\medskip

Then $-q_{n}<r_{n}<q_{n}$ and by Lemma \ref{l:31}, we have for every
$n\geq n_{1}$, $~c_{1} \in \Delta^{(n-1)}_{r_{n}}(c)$ and
$$\lim_{n\to+\infty}~b^{\prime}_{n}=0$$ In particular, for every $n\geq
n_{1}$,~ $c_{1}~~ \textrm{is an interior point of}~~
\Delta^{(n-1)}_{r_{n}}(c)\cap \Delta^{(n)}_{r_{n+1}}(c)$ and
$$-q_{n+1}<r_{n+1}-r_{n}<q_{n+1}.$$

- If $r_{n+1}-r_{n} \geq 0$, then by Lemma \ref{l:4611}  $$r_{n+1}-r_{n}\in
\{ q_{n-1}+jq_{n}: 0\leq j< a_{n+1}\}.$$

-If $r_{n+1}-r_{n}\leq 0$, then by Lemma \ref{l:4612}  $$r_{n+1}-r_{n}\in
\{ -jq_{n}: 1\leq j\leq a_{n+1}\}.$$

We distinguish four sub-cases.
\medskip

Case 2.1: $r_{n+1}-r_{n} \in \{q_{n-1}+jq_{n}:~0 \leq j <
a_{n+1}-1\}$. Then by the Denjoy's inequality, we have
$$l^{\prime}_{n} \leq e^{V}
 \frac{m(f^{-q_{n}}(\Delta^{(n)}_{r_{n}}(c)))}{m(\Delta^{(n-1)}_{r_{n}}(c))} \leq
 e^{V} b^{\prime}_{n}$$
(since $f^{-q_{n}}(\Delta^{(n)}_{r_{n}}(c))\subset
[f^{r_{n}}(c),~c_{1}]$).\\
\medskip

Case 2.2: $r_{n+1}-r_{n} = q_{n+1}-q_{n}$ (i.e. $j=a_{n+1}-1$).
 Then by the Denjoy's inequality, we have
\begin{align*}
l^{\prime}_{n} & =
 \frac{m(f^{-q_{n+1}+q_{n}}(\Delta^{(n)}_{r_{n+1}}(c)))}{m(\Delta^{(n-1)}_{r_{n}}(c))}\\
 & \leq
 e^{2V} \frac{m(\Delta^{(n)}_{r_{n+1}}(c))}{m(\Delta^{(n-1)}_{r_{n}}(c))} \\
 & \leq e^{2V}\left(\frac{m([f^{r_{n}}(c),~c_{1}])}{m(\Delta^{(n-1)}_{r_{n}}(c))}+
\frac{m([c_{1},~f^{r_{n+1}}(c)])}{m(\Delta^{(n)}_{r_{n+1}}(c))}\right)\\
& = e^{2V}(b^{\prime}_{n}+b^{\prime}_{n+1})
\end{align*}
(since $c_{1} \in \Delta^{(n)}_{r_{n+1}}(c) \subset
\Delta^{(n-1)}_{r_{n}}(c)$).
\bigskip

Case 2.3: $r_{n+1}-r_{n}\in \{-jq_{n}:~2\leq j\leq a_{n+1}\}$. Then
by the Denjoy's inequality, we have

\begin{align*}l^{\prime}_{n} & \leq e^{V}
 \frac{m(f^{-q_{n}}(\Delta^{(n)}_{r_{n}}(c)))}{m(\Delta^{(n-1)}_{r_{n}}(c))} \\
 & \leq
 e^{V} \frac{m([f^{r_{n}}(c),~c_{1}])}{m(\Delta^{(n-1)}_{r_{n}}(c))} \\
 & = e^{V}b^{\prime}_{n}
 \end{align*}
(since $f^{-q_{n}}(\Delta^{(n)}_{r_{n}}(c))\subset [f^{r_{n}}(c),~c_{1}]$).\\

Case 2.4: $r_{n+1}-r_{n}=-q_{n}$ (i.e. $j=1$). Then by the Denjoy's inequality,
we have
\begin{align*}
l^{\prime}_{n} & =
 \frac{m(f^{-q_{n}}(\Delta^{(n)}_{r_{n+1}}(c)))}{m(\Delta^{(n-1)}_{r_{n}}(c))}\\
 & \leq
 e^{V} \frac{m(\Delta^{(n)}_{r_{n+1}}(c))}{m(\Delta^{(n-1)}_{r_{n}}(c))} \\
 & \leq e^{V}\left(\frac{m([f^{r_{n}}(c),~c_{1}])}{m(\Delta^{(n-1)}_{r_{n}}(c))}+
\frac{m([c_{1},~f^{r_{n+1}}(c)])}{m(\Delta^{(n)}_{r_{n+1}}(c))}\right)\\
& = e^{V}(b^{\prime}_{n}+b^{\prime}_{n+1})
\end{align*}
(since $c_{1} \in \Delta^{(n)}_{r_{n+1}}(c) \subset
\Delta^{(n-1)}_{r_{n}}(c)$).
\bigskip

In either of the four cases [2.1$-$2.4] above, we conclude that for every $n\geq n_{0}$,~$$l^{\prime}_{n}\leq
e^{2V}(b^{\prime}_{n} + b^{\prime}_{n+1}).$$ So as in Case 1, we
conclude that  $\underset{n\to
+\infty}\lim ~a_{n+1} = +\infty$, a contradiction with that $\rho(f)$ is of bounded type.\\
\medskip

\textbf{Case 3}: $T(c_{1}, \mathbb{N}_{0})=\{0,1\}$. In this case,
there exists an infinite subset
$M_{1}$ of $\mathbb{N}_{0}$ such that $c_{1}\in Q_{2}(M_{1})\cap Q_{3}(M_{2})$, where $M_{2} = \mathbb{N}_{0}\backslash M_{1}$.
Hence there exists $n_{1}\in \mathbb{N}_{0}$ such that for all $n\geq
n_{1}, n\in M_{1}$, $c_{1}$ is an interior point
of $\Delta^{(n-1)}_{s_{n}}(c)$ and there
exists $n_{2}\in \mathbb{N}_{0}$ such that for all $n\geq
n_{2}, n\in M_{2}$, $c_{1}$ is an interior point of $\Delta^{(n-1)}_{r_{n}}(c)$. We distinguish the following cases.
\bigskip

\textit{Case 3.1}:  $(n, n+1)\in M_{1}\times M_{1}$ (resp. $(n, n+1)\in M_{2}\times M_{2}$) holds for infinitely many $n$, say
$n\in M_{1}^{\prime}\subset M_{1}$ (resp. $M_{2}^{\prime}\subset M_{2}$) with  $M_{1}^{\prime}$ (resp. $M_{2}^{\prime}$) infinite.
In this case, $c_{1}$ is an interior point of $\Delta^{(n-1)}_{s_{n}}(c)\cap \Delta^{(n)}_{s_{n+1}}(c)$ (resp.
$\Delta^{(n-1)}_{r_{n}}(c)\cap
\Delta^{(n)}_{r_{n+1}}(c)$) for $n\in M_{1}^{\prime}$ (resp. $n\in M_{2}^{\prime}$). Hence as in Cases 1 and 2,
we obtain that $$\underset{n\to
+\infty, n\in M_{1}^{\prime}}\lim ~a_{n+1} = +\infty \ (\textrm{resp}. \underset{n\to
+\infty, n\in M_{2}^{\prime}}\lim ~a_{n+1} = +\infty).$$
\bigskip

\textit{Case 3.2}: $(n, n+1)\in M_{1}\times M_{2}$ holds for infinitely many $n$. In this case, $c_{1}$ is an interior point of
$\Delta^{(n-1)}_{s_{n}}(c)\cap \Delta^{(n)}_{r_{n+1}}(c)$. Hence we have one of the following sub-cases.
\bigskip

 $A_{1}$: $r_{n+1}-s_{n}\in \{q_{n-1}+jq_{n}: ~1\leq j\leq a_{n+1}\}$. In this case, we have, by the Denjoy's inequality:
\begin{equation*}
l_{n} \leq e^{V}\frac{m([c_{1},~f^{s_{n}+q_{n-1}}(c)])}{m(\Delta^{(n-1)}_{s_{n}}(c))}=e^{V}b_{n}
\end{equation*}
(since $f^{q_{n-1}}(\Delta^{(n)}_{s_{n}}(c))\subset [c_{1},~f^{s_{n}+q_{n-1}}(c)])$.\\
\medskip

$A_{2}$:  $r_{n+1}-s_{n}\in \{-jq_{n}:~1\leq j\leq a_{n+1}-1\}$. In this case, \begin{equation*}
l_{n} \leq e^{3V}\frac{m([c_{1},~f^{s_{n}+q_{n-1}}(c)])}{m(\Delta^{(n-1)}_{s_{n}}(c))}=e^{3V}b_{n}
\end{equation*}
(since $f^{-(q_{n+1}-q_{n-1})}(\Delta^{(n)}_{s_{n}}(c))\subset
[c_{1},~f^{s_{n}+q_{n-1}}(c)])$.
\medskip

$A_{3}$:   $r_{n+1}-s_{n}= q_{n-1}$.
\bigskip

\textit{Case 3.3}:  $(n, n+1)\in M_{2}\times M_{1}$ holds for infinitely many $n$. In this case,
$c_{1}$ is an interior point of $\Delta^{(n-1)}_{r_{n}}(c)\cap \Delta^{(n)}_{s_{n+1}}(c)$. Hence we have one of the following sub-cases.
\bigskip

$B_{1}$: $s_{n+1}-r_{n}\in \{q_{n-1}+jq_{n}: ~0\leq j\leq a_{n+1}-2\}$. In this case, by the Denjoy's inequality, we have
\[l^{\prime}_{n} \leq e^{3V}b^{\prime}_{n} \ ( \textrm{since } \ f^{q_{n+1}-q_{n}}(\Delta^{(n)}_{r_{n}}(c))\subset
[f^{r_{n}}(c),~c_{1}]).\]
\medskip

$B_{2}$:  $s_{n+1}-r_{n}\in \{-jq_{n}: ~2\leq j\leq a_{n+1}+1\}$. In this case, by the Denjoy's inequality, we have
\[l^{\prime}_{n} \leq e^{V}b^{\prime}_{n} \ (\textrm{since} f^{-q_{n}}(\Delta^{(n)}_{r_{n}}(c))\subset [f^{r_{n}}(c),~c_{1}]).\]
\medskip

$B_{3}$:  $s_{n+1}-r_{n} = -q_{n}$.
\bigskip
\medskip

\textit{Claim.} One of the cases $A_{1}$, $A_{2}$, $B_{1}$ and $B_{2}$ holds for infinitely many $n$.
\medskip

Indeed, on the contrary, the cases $A_{3}$ and $B_{3}$ hold in particular for
infinitely many $n \in \mathbb{N}_{0}$. So for every $n\in M_{1}$, we
have $n+1\in M_{2}$ and $r_{n+1} = s_{n} + q_{n-1}$ and for every
$n\in M_{2}$, we have $n+1\in M_{1}$ and $s_{n+1}-r_{n}= -q_{n}$. So
for every $n\in M_{1}$, we have $n+1\in M_{2}$, $n+2\in M_{1}$ and
$n+3\in M_{2}$. Hence $r_{n+1}-s_{n+2}= q_{n+1}$ and
$r_{n+3}-s_{n+2}= q_{n+1}$, so $r_{n+1} = r_{n+3}$, for every $n\in
M_{1}$. Thus $(r_{n})_{n\in M_{2}}$ is a constant $r\in \mathbb{Z}$.
It follows that $c_{1}\in \Delta^{(n-1)}_{r}(c)$ for every $n\in
M_{2}$ and hence $m\left([f^{-r}(c_{1}), c]\right)\leq
m\left(\Delta^{(n-1)}_{0}(c)\right)$. Since $\underset{n\to
+\infty}\lim m\left(\Delta^{(n-1)}_{0}(c)\right) = 0$, so $c_{1} =
f^{r}(c)$, this contradicts the hypothesis $c_{1}\notin O_{f}(c)$.
\bigskip

If $A_{1}$ (resp. $A_{2}$) holds for infinitely many $n$, say $n\in M_{1}^{\prime\prime}\subset M_{1}$ (resp.
$M_{1}^{\prime\prime\prime}$)
with  $M_{1}^{\prime\prime}$ (resp. $M_{1}^{\prime\prime\prime}$) infinite,
so as in Cases 1 and 2,
we obtain that $$\underset{n\to
+\infty, n\in M_{1}^{\prime\prime}}\lim ~a_{n+1} = +\infty ~(\textrm{resp}. \underset{n\to
+\infty, n\in M_{1}^{\prime\prime\prime}}\lim ~a_{n+1} = +\infty).$$
\medskip
\medskip
If $B_{1}$ (resp. $B_{2}$) holds for infinitely many $n$, say $n\in M_{2}^{\prime\prime}\subset M_{2}$
(resp. $M_{2}^{\prime\prime\prime}\subset M_{2}$) with
$M_{2}^{\prime\prime}$ (resp. $M_{2}^{\prime\prime\prime}$) infinite, so similarily as for $A_{1}$ and  $A_{2}$,
we obtain that $$\underset{n\to
+\infty, n\in M_{2}^{\prime\prime}}\lim ~a_{n+1} = +\infty ~ (\textrm{resp}. \underset{n\to
+\infty, n\in M_{2}^{\prime\prime\prime}}\lim ~a_{n+1} = +\infty.)$$\\

In either of the three cases [3.1$-$3.3] above, we get a contradiction with that
$\rho(f)$ is of bounded type. We conclude that the Case 3 doesn't occur and so $Q_{4}( \mathbb{N}_{0})=\emptyset$. \qed
%
%

\bigskip
\bigskip

In the sequel, 
we write for simplicity $\mathbb{N}_{0}$ instead of $\mathbb{N}_{c,c_{1}}$ and we denote by

 $$E(c_{1},c):= \{c\}\cup Q_{2}(\mathbb{N}_{0})\cup Q_{3}(\mathbb{N}_{0})$$

So one has
$$B\backslash\{c\} = Q_{1}( \mathbb{N}_{0})\ \amalg~Q_{2}( \mathbb{N}_{0}) \
\amalg~Q_{3}(\mathbb{N}_{0})  .$$
\medskip

For $a \in S^{1},~n \in \mathbb{N}_{0}$ and $\gamma >0,$ set
\medskip

\textbullet ~ $~l_{n}: = m(\Delta^{(n)}(y_{1}))$\\

\textbullet ~ $V^{(f)}_{n,\gamma}(a):=[a-\gamma ~l_{n-1},~a+\gamma
~l_{n-1}]$ \\

 \textbullet ~ $k_{n}(d)=\begin{cases}

                      i_{n}(d), & \textrm{if}~d \in \{c\} \cup Q_{2}(\mathbb{N}_{0}) \\
                      j_{n}(d), & \textrm{if}~d \in Q_{3}(\mathbb{N}_{0})
                                               \end{cases}$
\medskip
\medskip

\begin{prop}\label{p:47} There exists a positive constant $\gamma_{0}>0$ such that:
\medskip

\begin{enumerate}
  \item For every $d\in E(c_{1},c)$ and every $0<\gamma< \gamma_{0}$, there exists
$n_{\gamma} \in \mathbb{N}_{0}
  $ such that for any $n \in \mathbb{N}_{0},~n \geq n_{\gamma}$ there exists a unique integer $0 \leq k_{n}(d)<q_{n}$ such that

  $$f^{-k_{n}(d)}(d) \in V^{(f)}_{n,\gamma}(f^{-i_{n}(c)}(c))\subset  ]y_{1},y_{3}[$$

  \item For every $n \in \mathbb{N}_{0},~\bigcup\limits_{i=0}^{q_{n}-1} f^{-i}\big(B\backslash E(c_{1},c)\big)\subset S^{1}\backslash
  V^{(f)}_{n,\gamma_{0}}(f^{-i_{n}(c)}(c))$.\\
\item For every $~d \in Q_{11}( \mathbb{N}_{0})$ \textrm{and } $n\in \mathbb{N}_{0}$, \textrm{we have } $f^{-i_{n}(d)}(d)\in
~]f^{q_{n}}(y_{1})+\gamma_{0}l_{n-1}, y_{2}-\gamma_{0}l_{n-1}[$.\\
   \item For every $~d \in Q_{12}( \mathbb{N}_{0})$ \textrm{and }$n\in \mathbb{N}_{0}$, \textrm{we have }
$f^{-j_{n}(d)}(d) \in ~]y_{2}+\gamma_{0}l_{n-1},
y_{3}-\gamma_{0}l_{n-1}[$.
\end{enumerate}
\end{prop}
\medskip

\begin{proof} Let $d\in B\backslash \{c\}$.

Assertion (1): The
intervals $f^{i_{n}(d)}([y_{1},~y_{2}])$ and
$f^{j_{n}(d)}([y_{2},~y_{3}])$ are comparable (that is there is a
positive constant $0<K<1$ such that $K \leq
\dfrac{m(f^{i_{n}(d)}([y_{1},~y_{2}]))}{m(f^{j_{n}(d)}([y_{2},~y_{3}]))}
\leq K^{-1}$ for every $n \in \mathbb{N}_{0}$) and are comparable to
their union $[f^{j_{n}(d)}(y_{2}),~f^{i_{n}(d)}(y_{2})]$ which
implies, by Lemma \ref{l:31} and Denjoy's inequality, that the
four intervals
$[y_{1},~y_{2}],~[y_{2},~y_{3}],~[f^{j_{n}(d)-i_{n}(d)}(y_{2}),y_{2}],~[y_{2},~f^{i_{n}(d)-j_{n}(d)}(y_{2})]$
are comparable. In particular, the ratios
$\dfrac{m([f^{-i_{n}(d)}(d),~y_{2}])}{m([y_{1},~y_{2}])}$ and
$\dfrac{m([f^{-i_{n}(d)}(d),~y_{2}])}{m([f^{j_{n}(d)-i_{n}(d)}(y_{2}),~y_{2}])}=t_{n}(d)$ are comparable.
So, for $d\in Q_{2}(\mathbb{N}_{0}),~\lim_{n \to +\infty}
~t_{n}(d)=0$ which is equivalent to the assertion:
$$\forall ~\gamma >0,~\exists ~ n_{\gamma} \in \mathbb{N}_{0}
~\textrm{such that for every }~n \geq n_{\gamma},~ n \in
\mathbb{N}_{0}~:~f^{-i_{n}(d)}(d) \in V^{(f)}_{n,\gamma}(y_{2}).
$$
$\bullet$ ~ Similarly, the ratios
$\dfrac{m([y_{2},~f^{-j_{n}(d)}(d)])}{m([y_{1},~y_{2}])}$ and
$\dfrac{m([y_{2},~f^{-j_{n}(d)}(d)])}{m([y_{2},~f^{i_{n}(d)-j_{n}(d)}(y_{2})])}$
are comparable, which is comparable to the
ratio
$\dfrac{m([f^{j_{n}(d)-i_{n}(d)}(y_{2}),~f^{-i_{n}(d)}(d)])}{m([f^{j_{n}(d)-i_{n}(d)}(y_{2}),~y_{2}])}
=1-t_{n}(d)
$ (by the Lemma \ref{l:31}). So, for $d \in Q_{3}(\mathbb{N}_{0}),~\lim_{n \to +\infty}
~t_{n}(d)=1$, which is equivalent to the assertion:
$$\forall ~\gamma >0,~\exists ~ n_{\gamma}\in \mathbb{N}_{0}
~\textrm{such that for every }~n \geq n_{\gamma},~ n \in
\mathbb{N}_{0}:~f^{-j_{n}(d)}(d) \in V^{(f)}_{n,\gamma}(y_{2}).
$$

Assertion (2): Let $d \in Q_{1}(\mathbb{N}_{0})$. There exists
$0<\gamma (d)<1$ such that for every $n \in \mathbb{N}_{0},~t_{n}(d)
\geq \gamma (d) ~ \textrm{and}~1-t_{n}(d) \geq  \gamma (d) $. Since
$Q_{1}(\mathbb{N}_{0})$ is finite there exists $0<\gamma_{0}<1$ and
an infinite subset $M_{2}$  of $\mathbb{N}_{0}$
($\mathbb{N}_{0}\backslash M_{2}$ is finite) such that the points
$f^{-i_{n}(d)}(d)$ and $f^{-j_{n}(d)}(d)$ are contained in the set
$S^{1}\backslash V^{(f)}_{n,\gamma_{0}}(y_{2})$ for all $d\in
Q_{1}(M_{2})$  and $n \in M_{2}$. By Lemma \ref{l:25}, $~ d\notin f^{i}
([y_{1},~y_{2}])\cup f^{j} ([y_{2},~y_{3}])$, for every $0
\leq i < q_{n}, ~i \neq i_{n}(d), ~j = \varphi_{n}(i)$. Hence, for every $0
\leq i <q_{n}, ~i \neq i_{n}(d):~f^{-i}(d) \notin [y_{1},~y_{2}]$
and every $0 \leq j <q_{n},~j \neq
j_{n}(d)=\varphi_{n}(i_{n}(d)): f^{-j}(d)\notin [y_{2},~y_{3}]$. In
particular, for every $0 \leq i <q_{n},~i \neq i_{n}(d),~i \neq
j_{n}(d): ~f^{-i}(d) \notin [y_{1},~y_{3}]$. So $f^{-i}(d) \notin
V^{(f)}_{n,\gamma_{0}}(y_{2})$ since
$V^{(f)}_{n,\gamma_{0}}(y_{2})\subset [y_{1},~y_{3}]$. Then
$\bigcup_{0 \leq i <q_{n}}f^{-i}(B\backslash E(c_{1},c))\subset
S^{1}
\backslash V^{(f)}_{n,\gamma_{0}}(y_{2})$.
\bigskip

Write $\mathbb{N}_{0}$ instead of $M_{2}$.
\
\\

Assertion (3): Let $d \in Q_{11}(\mathbb{N}_{0})$. Then there exist
an integer ~ $n(d) \in \mathbb{N}_{0}$ such that for every $n\in
\mathbb{N}_{0}~( n \geq n(d))~,~0 \leq i_{n}(d)< q_{n-1}$,

\begin{align*}
f^{-i_{n}(d)}(d) \in [y_{1},~y_{2}]\cup
f^{j_{n}(d)-i_{n}(d)}([y_{2},~y_{3}])&=[y_{1},~y_{2}]\cup
[f^{q_{n}}(y_{1}),~f^{q_{n}}(y_{2})]\\\\
& =[f^{q_{n}}(y_{1}),y_{2}]
\end{align*}
 and
$$\gamma_{0} \leq
t_{n}(d)=\dfrac{m([f^{-i_{n}(d)}(d),~y_{2}])}{m([f^{q_{n}}(y_{1}),~y_{2}])}
\leq 1-\gamma_{0}.$$
\medskip
\
\\
Since the intervals $[f^{q_{n}}(y_{1}),~y_{2}]$ and $[y_{1},~y_{2}]$
are comparable, there exists ~ $0 <\gamma_{1}<\gamma_{0}$ and an
integer $n^{\prime}(d) \in \mathbb{N}_{0}$ such that for all \ $d\in
\mathbb{N}_{0}$ and $n \geq n^{\prime}(d)$: $f^{-i_{n}(d)}(d) \in
[f^{q_{n}}(y_{1})+\gamma_{1}l_{n-1},~y_{2}-\gamma_{1}l_{n-1}]$.
Since $Q_{11}(\mathbb{N}_{0})$ is finite, there exists an integer
$n_{1} \in \mathbb{N}_{0}$ such that for all \ $d\in Q_{11}(
\mathbb{N}_{0})$ and $n\in \mathbb{N}_{0},~n \geq n_{1}$:
$$f^{-i_{n}(d)}(d) \in
[f^{q_{n}}(y_{1})+\gamma_{1}l_{n-1},~y_{2}-\gamma_{1}l_{n-1}].$$
\medskip
\
\\
Assertion (4): Let $d\in Q_{12}( \mathbb{N}_{0})$. Then there exist
an integer ~ $n(d) \in \mathbb{N}_{0}$ such that for every $n\in \mathbb{N}_{0}, \ n \geq n(d)$,  $q_{n-1} \leq i_{n}(d) <
q_{n}$, $f^{-j_{n}(d)}(d) \in [y_{2},~y_{3}]$ and
$$\gamma_{0} \leq
t_{n}(d) = \dfrac{m([f^{-i_{n}(d)}(d),~y_{2}])}{m([y_{1},~y_{2}])}
\leq 1-\gamma_{0}.$$
\medskip
\
\\
Since the ratios
$$\dfrac{m([f^{-j_{n}(d)}(d), ~y_{3}])}{m([y_{2},~y_{3}])}
~~~~~~~~~ \  \textrm{ and } \ ~~~~~~~~~~~
\dfrac{m([f^{-i_{n}(d)}(d),~y_{2}])}{m([y_{1},~y_{2}])}$$
\medskip
\
\\
 are comparable, there exists $0 <\gamma_{2}<\gamma_{0}$ such that
$f^{-j_{n}(d)}(d) \in
[y_{2}+\gamma_{2}l_{n-1},~y_{3}-\gamma_{2}l_{n-1}], \ \textrm{for
all } \ d\in Q_{12}(\mathbb{N}_{0})$. Since $Q_{12}(\mathbb{N}_{0})$
is finite, there exists an integer ~ $n_{2} \in \mathbb{N}_{0}$ such
that $\textrm{for all } \ d\in Q_{12}(\mathbb{N}_{0})$ and $n\in
\mathbb{N}_{0},~n \geq n_{2}$:  $$f^{-j_{n}(d)}(d) \in
[y_{2}+\gamma_{2}l_{n-1},~y_{3}-\gamma_{2}l_{n-1}]$$
\end{proof}

\begin{cor}\label{c:37}
For every $d \in E(c_{1},c)$,  we have \[\lim_{n \to + \infty,~n \in
\mathbb{N}_{0}}\dfrac{m([f^{-k_{n}(d)}(d),~y_{2}])}{m([y_{1},~y_{2}])}
=0 \]
\end{cor}
\medskip

\begin{proof} It is a consequence of Proposition \ref{p:47}.
\end{proof}
\medskip

\
\\
\section{\bf Secondary Cells}

We introduce in this section a new triple $(z_{1},z_{2},z_{3})$ in
$U_{\delta}(x_{0})$ depending on two real parameters $\beta$ and
$\gamma$. Let $\beta,~\gamma \in ]0,\gamma_{0}[~(\beta <\gamma)$.
For $n \in \mathbb{N}_{0}~(n \geq n_{\beta}),~$ we define the points
$z_{1},z_{2},z_{3}\in S^{1}$ ($z_{1}\prec z_{2}\prec z_{3})$ as
follows: $z_{2} = y_{2}$ and $z_{1}, ~z_{3}\in
V^{(f)}_{n,\gamma}(z_{2})\backslash V^{(f)}_{n,\beta}(z_{2})$ i.e.

$$\beta\leq \frac{m([z_{1},z_{2}])}{m([y_{1},y_{2}])}\leq
\gamma ~\textrm{and }~\beta\leq
\frac{m([z_{2},z_{3}])}{m([y_{1},y_{2}])}\leq \gamma $$
\bigskip

\textbf{Definition (secondary cell)}. The triple
$(z_{1},~z_{2},~z_{3})$ is called a $(\beta,\gamma)$-derived cell
associated to $(f,x_{0},c,\mathbb{N}_{0},\delta)$.
\medskip

\begin{prop}\label{p:42} Under the notations above, for every $n \in \mathbb{N}_{0},~n \geq n_{\beta} $ we have:
\begin{itemize}
 \item [(0)] $[z_{1},z_{3}]\subset  V^{(f)}_{n,\gamma}(y_{2})\subset [y_{1},y_{3}]$.

 \item [(1)] For any $d\in E(c_{1},c)$, there exists a unique $0\leq k_{n}(d)<q_{n}$ such that $f^{-k_{n}(d)}(d)\in  V^{(f)}_{n,\frac{1}{4}\beta}(y_{2})$.

 \item [(2)] $~m\left(f^{j}([z_{1},z_{3}])\right)\leq K\lambda^{n}$, for every $0\leq j<q_{n}$, where
$\lambda = (1+e^{-V})^{-\frac{1}{2}}$ and  $K$ is a constant
independent of $n$, $V= \textrm{Var}(\log(\mathrm{Df}))$.

 \item [(3)] $~\sum_{j=0}^{q_{n}-1}m\left(f^{j}([z_{1},z_{3}])\right)\leq
2$.

 \item [(4)] The triples  $\big(z_{1},z_{2},z_{3}\big)$~ and
$(f^{q_{n}}\big(z_{1}),f^{q_{n}}(z_{2}),f^{q_{n}}(z_{3})\big)$
satisfy conditions  $(a)$  and  $(b)$ of Definition \ref{d:32} for
$x_{0}$  and the constant $R(\beta)  = ~\frac{Re^{2V}}{\beta}$,
where $R=e^{3V}+e^{V}+1$.

 \item [(5)] $f^{i}([z_{1},z_{3}])$ does not contain any  point of $B$, for every $0\leq i< q_{n}$ and $i\neq k_{n}(d),~d \in
E(c_{1},c)$.

\item [(6)] For every $d\in E(c_{1},c),~d \in
f^{k_{n}(d)}(]z_{1},z_{3}[)$.

\end{itemize}
\end{prop}
\medskip

\begin{proof}
Assertion (0): This results from the definition of $\big(z_{1},z_{2},z_{3}\big)$, for $\gamma$ small.\\

Assertion (1): This follows from the Assertion (1) in Proposition \ref{p:47}.\\

Assertions (2) and (3) hold since $[z_{1},z_{3}]\subset
[y_{1},y_{3}]$ and by applying the same proof to
$\big(y_{1},y_{2},y_{3}\big)$ instead of
$\big(y_{1},y_{2},y_{3}\big)$ of the properties (c-2) and (c-3) of Proposition \ref{p:41}.\\
Assertion (4): We have
$$\dfrac{\beta}{\gamma} = \dfrac{\beta m([y_{1},y_{2}])}{\gamma
m([y_{1},y_{2}])} \leq \dfrac{m([z_{2},z_{3}])}{m([z_{1},z_{2}])}
\leq \dfrac{\gamma m([y_{1},y_{2}])}{\beta m([y_{1},y_{2}])} =
\dfrac{\gamma}{\beta}.$$

Then by Lemma \ref{l:31}, we get
\bigskip

\begin{align*}
  e^{-2V}\dfrac{\beta}{\gamma} & = e^{-2V}\dfrac{\beta m([y_{1},y_{2}])}{\gamma m([y_{1},y_{2}])} \\
& \leq e^{-2V}\dfrac{m([z_{2},z_{3}])}{m([z_{1},z_{2}])} \leq
\dfrac{m(f^{q_{n}}([z_{2},z_{3}]))}{m(f^{q_{n}}([z_{1},z_{2}]))}\\
& \leq  e^{2V}\dfrac{m([z_{2},z_{3}])}{m([z_{1},z_{2}])} \\
& \leq e^{2V}\dfrac{\gamma m([y_{1},y_{2}])}{\beta m([y_{1},y_{2}])}
= e^{2V} \dfrac{\gamma}{\beta}.
\end{align*}
\\\\

 Hence $\big(z_{1},z_{2},z_{3}\big)$ satisfies conditions (a) and (b) of Definition \ref{d:32} for the
constant $R(\beta, \gamma) = \dfrac{R}{\beta}e^{2V}$. On the other
hand, since $[z_{1},z_{3}]\subset [f^{q_{n}}(y_{1}),y_{3}],~x_{0}\in
[f^{q_{n}}(y_{1}),y_{3}]$ and

$\max \left(m([f^{q_{n}}(y_{1}),x_{0}]); m([x_{0},y_{3}])\right)\leq R
m([y_{1},y_{2}]),$ it follows that $$\underset{1\leq i\leq 3}\max(
m([x_{0},z_{i}])\leq
\max\left(m([f^{q_{n}}(y_{1}),x_{0}]);m([x_{0},y_{3}])\right)\leq R
m([y_{1},y_{2}])\leq \dfrac{R}{\beta}m([z_{1},z_{2}]).$$

Furthermore, as $[z_{1},z_{3}]\subset ] y_{1},y_{3}[ $ then
$f^{q_{n}}([z_{1},z_{3}])\subset f^{q_{n}}(] y_{1},y_{3}[)$
\\

Since $x_{0}\in [f^{q_{n}}(y_{1}),y_{3}]~ \textrm{and} ~ f^{q_{n}}(]
y_{1},y_{3}[)\cup f^{q_{n-1}} ( ]f^{q_{n}}(y_{1}),y_{2}[ )\cup
f^{q_{n}-q_{n-1}}( ] y_{2},y_{3}[)\subset [f^{q_{n}}(y_{1}),y_{3}],$

we get

\begin{align*}
\underset{1\leq i\leq 3}\max m([f^{q_{n}}(z_{i}),x_{0}])&\leq
\max\left (m([f^{q_{n}}(y_{1}),x_{0}]);m([x_{0},y_{3}])\right )\\
& \leq
\dfrac{R}{\beta}m([z_{1},z_{2}]) \\
& \leq \dfrac{R}{\beta}e^{2V}m(f^{q_{n}}([z_{1},z_{2}])~~
(\textrm{by Denjoy's inequality} ).
\end{align*}
 Therefore

$\big(f^{q_{n}}(z_{1}),f^{q_{n}}(z_{2}),f^{q_{n}}(z_{3})\big)$
satisfies conditions (a) and (b) of Definition \ref{d:32} for the
constant $R(\beta,\gamma) = \dfrac{R}{\beta}e^{2V}$. \\

Assertion (5): If $d\in  B$ and $0\leq i < q_{n}$ be an integer such
that $d\in f^{i}([z_{1},z_{3}])$ then $f^{-i}(d)\in
[z_{1},z_{3}]\subset V^{(f)}_{n,\gamma}(y_{2})\subset
V^{(f)}_{n,\gamma_{0}}(y_{2})$. According to the Assertion (2) of
Propositions \ref{p:47}, we obtain $f^{-i}(d)\notin
\bigcup\limits_{j=0}^{q_{n}-1} f^{-j}\big(B\backslash E(c_{1},c)
\big)$, $0\leq i < q_{n}$. In particular, $f^{-i}(d)\notin
f^{-i}\big(B\backslash E(c_{1},c) \big)$. Hence $d\in E(c_{1},c)$
such that $0\leq i < q_{n}$ and $f^{-i}(d) \in
V^{(f)}_{n,\gamma}(y_{2})$. By Assertion (1), $i=k_{n}(d)$.\\

Assertion (6) is a consequence of Assertions (0) and (1).
\end{proof}
\medskip

\begin{prop}[\cite{A}, Proposition 4.2]\label{p:46} Let $\left(y_{1},y_{2},y_{3}\right)\subset
U_{\delta}(x_{0})$ be the primary cell and $(z_{1},z_{2},z_{3})$ an
$(\beta,\gamma)$ secondary cell associated to
$(f,x_{0},c,\mathbb{N}_{0})$. Let $h$ be the homeomorphism
conjugating $f$ to $g$ and assume that it admits in $x_{0}$ a positive
derivative $\textrm{Dh}(x_{0})>0$. Then the following properties hold:

\begin{itemize}
 \item [(1)] $(h(y_{1}),h(y_{2}),h(y_{3}))$ is a
primary cell associated to $(g,h(x_{0}),h(c),\mathbb{N}_{0})$.

\item [(2)] $(h(z_{1}),h(z_{2}),h(z_{3}))$ is a $(\frac{\beta}{2},2\gamma)$-derived cell associated to
$(g,h(x_{0}),h(c),\mathbb{N}_{0})$.

\item [(3)] There is an integer $n_{\gamma} \in \mathbb{N}_{0}$ such
that for every $n \in \mathbb{N}_{0},~n \geq n_{\gamma}$ :
$$h \big(V^{(f)}_{n,\gamma}(y_{2}) \big) \subset V^{(g)}_{n,2\gamma}\big(h(y_{2})\big)
\subset h \big(V^{(f)}_{n,4\gamma}(y_{2}) \big)$$
\end{itemize}
\end{prop}
\medskip

\begin{prop}\label{p:43} Under the notations of Proposition \ref{p:46}, for every $n
\in \mathbb{N}_{0},~n \geq n_{\beta} $ we have:
\medskip
\begin{itemize}

 \item [(h-0)] $[h(z_{1}),h(z_{3})]\subset  V^{(g)}_{n,2\gamma}(h(y_{2}))\subset
[h(y_{1}),h(y_{3})]$.

 \item [(h-1)] For any $d \in E(c_{1},c)$, $k_{n}(d)$ is the unique integer in $ [0,~q_{n}[$ such that
$g^{-k_{n}(d)}(h(d))\in  V^{(g)}_{n,\frac{1}{2}\beta}(h(y_{2}))$.

 \item [(h-2)] $~m\left(g^{j}([h(z_{1}),h(z_{3})])\right)\leq K^{\prime}(\lambda^{\prime})^{n}$, for every  $0\leq
j<q_{n}$,
where $\lambda^{\prime} = (1+e^{-V^{\prime}})^{-\frac{1}{2}}$ and
$K^{\prime}$ is a constant independent of  $n$,
$V^{\prime}=\textrm{Var}(\log(Dg))$.

 \item [(h-3)] $~\sum_{j=0}^{q_{n}-1}m\left(g^{j}([h(z_{1}),h(z_{3})])\right)\leq 2$

 \item [(h-4)] $g^{i}([h(z_{1}),h(z_{3})])$ does not contain any break point of $h(B)$ for every $0\leq i< q_{n}$ and $i\neq k_{n}(d),
~d \in E(c_{1},c)$.

\item [(h-5)] For every $d\in E(c_{1},c)$,

$$ \lim_{n\to
+\infty}\dfrac{m\big([g^{-k_{n}(d)}(h(d)),h(y_{2})]\big)}{m
\big([h(y_{1}),h(y_{2}) ] \big)} = 0 $$

\item [(h-6)] For every $d \in E(c_{1},c),~h(d) \in
g^{k_{n}(d)}(]h(z_{1}),h(z_{3})[)$.

\end{itemize}
\end{prop}
\medskip

\begin{proof} Assertion $($h-0$)$ is a consequence of Proposition
\ref{p:46}, (3).
\
\\

(h-1$)$: Let $d\in E(c_{1},c)$, by Proposition \ref{p:42} (1), there
exists  $0\leq k_{n}(d)<q_{n}$ such that $f^{-k_{n}(d)}(d)\in
V^{(f)}_{n,\frac{1}{4}\beta}(z_{2})$; then $g^{-k_{n}(d)}(h(d)) =
h(f^{-k_{n}(d)}(d))\in h(V^{(f)}_{n,\frac{1}{4}\beta}(z_{2}))\subset
V^{(g)}_{n,\frac{1}{2}\beta}(h(z_{2}))$ for every $n\geq
n^{\prime\prime}_{\beta},~~n \in \mathbb{N}_{0}$. Uniqueness: If
~~$0 \leq i,j < q_{n}$ satisfy $g^{-i}(h(d)),~g^{-j}(h(d))\in
V^{(g)}_{n,2\gamma}(h(z_{2}))$, then $f^{-i}(d), ~f^{-j}(d)\in
h^{-1}\big( V^{(g)}_{n,\gamma}(h(z_{2})\big)\subset V^{(f)}_{n,2\gamma}(z_{2})$. Therefore by Proposition \ref{p:42} (1), $i=j$.\\
\
\\
$($h-2$)$ and $($h-3$)$: by Proposition \ref{p:46}, (2),
$(h(z_{1}),h(z_{2}),h(z_{3}))$ is a
$(\frac{\beta}{2},2\gamma)$-derived cell associated to
$(g,h(x_{0}),h(c),\mathbb{N}_{0})$. Therefore, the proof is similar to Proposition \ref{p:42}, \textit{(2)} and \textit{(3)} by replacing $f$ by $g$.\\
\
\\
$($h-4$)$: Let $0\leq i<q_{n}$ and assume that
$g^{i}\big([h(z_{1}),h(z_{3})] \big)$ contains a break point
$d^{\prime} = h(d)$ of $g$ for a break point $d$ of $f$, then
$f^{i}\big([z_{1},z_{3}] \big)$ contains the break point $d$ of $f$.
By Proposition \ref{p:42}, (1) and (5), we have $d\in
E(c_{1},c)$ and $i = k_{n}(d)$.
\
\\
$($h-5$)$: Let  $d \in E(c_{1},c)$ and $\gamma >0$. By Proposition
\ref{p:41} (6), there is an integer $n_{\gamma}^{\prime} \in \mathbb{N}_{0}$
such that for every $n \in \mathbb{N}_{0},~n \geq n_{\gamma}^{\prime}$ :
$f^{-k_{n}(d)}(d) \in V^{(f)}(y_{2})$. On the other hand, by
Proposition \ref{p:46}, (3), there is an integer $n_{\gamma} \in
\mathbb{N}_{0},~n_{\gamma}\geq n_{\gamma}^{\prime}$ such that for every $n
\in \mathbb{N}_{0},~n \geq n_{\gamma}$ : $ h
\big(V^{(f)}_{n,\frac{\gamma}{2}}(y_{2}) \big) \subset
V^{(g)}_{n,\gamma}(h(y_{2}))$. As $g^{-k_{n}(d)}(h(d)) =
h(f^{-k_{n}(d)}(d))$, it follows that $g^{-k_{n}(d)}(h(d))\in
V^{(g)}_{n,\gamma}(h(y_{2}))$ for every $n \geq n_{\gamma},~n \in
\mathbb{N}_{0}$. Since $\gamma>0$ is arbitrary, so, $$ \lim_{n \to
+\infty} \dfrac{m\big([g^{-k_{n}(d)}(h(d)),h(y_{2}) ] \big)}{m
\big([h(y_{1}),h(y_{2}) ] \big)}=0
$$
$($h-6$)$ is a consequence of Assertions ($h-0$) and ($h-1$).
\end{proof}
\medskip

\section{\bf Control of distortions}
\medskip

\begin{prop}[\cite{A}, Proposition 5.1]\label{p:2} Let $\mathbb{N}_{0}$ and $\gamma_{0}$ are as in Proposition \ref{p:47}. Let
$\beta,~\gamma \in ]0,\gamma_{0}[~(\beta <\gamma)$. Assume that the
conjugation homeomorphism $h$ from $f$ to $g$ admits at $x_{0}$ a
positive derivative $\mathrm{Dh}(x_{0}) = \omega_{0} >0$. Then there
exists an integer $n_{\beta} \in \mathbb{N}_{0}$ such that: for
every $n\in \mathbb{N}_{0}, \ n\geq n_{\beta},$ we have

$$\left|
\dfrac{\mathrm{Dr}_{g^{q_{n}}}\big(h(z_{1}),h(z_{2}),h(z_{3})\big)}{\mathrm{Dr}_{f^{q_{n}}}(z_{1},z_{2},z_{3})}
- 1\right|\leq \dfrac{2}{1-\gamma_{0}}\beta$$
\end{prop}
\bigskip
\
\\

The next proposition is an another distortion control which is
opposite to the one of Proposition \ref{p:2}, this allows us to
prove that the conjugation from $f$ to $g$ is singular with respect
to the Lebesgue measure.

\begin{prop}
\label{p:44} Assume that the irrational rotation number
of $f$ is of bounded type. Let $c, c_{1}$ in $B(f)$,  $c\neq c_{1}$,  $\mathbb{N}_{0}$ and $\gamma_{0}$ are as in
Proposition \ref{p:47}. Then
for any $\varepsilon >0$, there exists $0<\gamma<\gamma_{0}$ such
that for any $\beta \in ]0,\gamma[$ there exists $n_{\beta} \in
\mathbb{N}_{0}$ such that for all $n \geq n_{\beta},~n \in
\mathbb{N}_{0}$ the ($\beta,\gamma$) secondary cell
$(z_{1},z_{2},z_{3})$ associated to
$(f,x_{0},c,\mathbb{N}_{0},\delta)$ satisfies the following
inequality
$$ \left| \dfrac{\textrm{Dr}_{g^{q_{n}}}\left(h(z_{1}),h(z_{2}),h(z_{3})\right)}{\textrm{Dr}_{f^{q_{n}}}(z_{1},z_{2},z_{3})}-
\Pi(c_{1},c) \right|\leq A \ \varepsilon, ~~\textrm{for some
constant}  \; A> 0,$$

where
$$\Pi(c_{1},c)=
 \underset{d\in E (c_{1},c)}\prod \dfrac{\sigma_{g}(h(d))}{\sigma_{f}(d)}.$$
\end{prop}
\medskip

The proof of the proposition \ref{p:44} is an elaboration of the proof of (\cite{A}, Proposition 5.3) and so we only
describe the changes that are necessary.
\medskip

\begin{lem}\label{l:49q}
Assume that the
rotation number $\alpha$ of $f$ is of bounded type. Let $c_{1} \in C(f)\backslash \{c\}$, $u_{0}$, $\gamma_{c,c_{1}}$ and
$\mathbb{N}_{0}$ as in Proposition \ref{p:47}. Let $h$ be the conjugating from $f$ to $R_{\alpha}$:
$h\circ f= R_{\alpha}\circ h$. Assume that Dh$(x_{0})>0$. Then for any $\varepsilon >0$,
there exists $\eta >0$ such that for any $u\in ]0,\eta[$ there exists $n_{u,r} \in
\mathbb{N}_{c,c_{1}}$ such that for any $n \in \mathbb{N}_{c,c_{1}},\ n\geq
n_{u,r}$: the $($r,u$)$-derived cell $(z_{1},z_{2},z_{3})$ associated to
$(f,x_{0},c)$ satisfies
\
\\
$$\arrowvert\textrm{Dcr}_{f^{q_{n}}}(z_{1},z_{2},z_{3})-
\Pi_{f}(c,c_{1})\arrowvert\leq C_{1}\varepsilon,$$ where $$\Pi_{f}(c,c_{1})= \underset{b\in E (c_{1},c)}\prod
\sigma_{f}(b)$$ and \; $C_{1}$ is a positive constant.
\end{lem}
\bigskip

\begin{proof}
Let $\varepsilon>0$. Since $\textrm{Df}$ is absolutely continuous on
every interval $[c_{i},c_{i+1}]~(0 \leq i \leq p)$, there exists a
real $0<\eta_{0}<e^{2V}\gamma_{0}$ such that for every disjoint
intervals $([a_{j},b_{j}])_{j \in J}$ of continuity intervals of
$\textrm{Df}$ (these later are the $[c_{i},c_{i+1}]~(0 \leq i \leq
p)$) such that if $\sum_{j \in J}|a_{j}-b_{j}|<\eta_{0}$ then
$\sum_{j \in J}|\textrm{Df}(a_{j})-\textrm{Df}(b_{j})|<\varepsilon$.
Let $0<\beta<\gamma<\frac{1}{5}e^{-2V}\eta_{0}$ and
$(z_{1},z_{2},z_{3})$ be the $(\beta , \gamma)$-derived cell
associated
to $(f,x_{0},c,\mathbb{N}_{0})$. For $n\in \mathbb{N}_{0}$, denote the following:\\

 \textbullet ~ $I:= \{0,1,\dots,q_{n}-1\}\backslash \{k_{n}(d):~d \in E(c_{1},c)\}$.\\

\textbullet ~ $D_{i}(f):=
\textrm{Dcr}_{f}\big(f^{i}(z_{1}),f^{i}(z_{2}),f^{i}(z_{3})\big),~0
\leq
i<q_{n}$.\\

\textbullet ~ $P:= \prod_{i\in I}D_{i}(f)~~;~~ P_{1}:= \prod_{d\in
E(c_{1},c)}D_{k_{n}(d)}(f)$.\\

\textbullet ~ $D_{n}(f):=
\textrm{Dcr}_{f^{q_{n}}}(z_{1},z_{2},z_{3})$.
\medskip

One has $D_{n}(f) = PP_{1}$.
\bigskip

\textbf{Step 1}. We consider $D_{i}(f),~i \in I$. By Proposition
\ref{p:42} (5), $\textrm{Df}$ is continuous on the interval
$[f^{i}(z_{1}),f^{i}(z_{3})]$. Then by the mean value theorem, there
exist $a_{i}\in ]f^{i}(z_{1}),f^{i}(z_{2})[$ and $b_{i}\in
]f^{i}(z_{2}),f^{i}(z_{3})[$ such that $D_{i}(f)=
\dfrac{\textrm{Df}(a_{i})}{\textrm{Df}(b_{i})} $. Since
$\big([a_{i},f^{i}(z_{2})] \big)_{i\in I}$ are disjoint intervals of
continuity of $\textrm{Df}$, by Lemma \ref{l:31}, they satisfy the
inequalities :\\\\

$\begin{array}{ll}
   \sum_{i\in I}m([a_{i},f^{i}(z_{2})]) & \leq \sum_{i \in I}m([f^{i}(z_{1}),f^{i}(z_{2})])
\\\\

   & \leq \sum_{i\in I} e^{2V} \dfrac{m([z_{1},z_{2}])}{m([y_{1},y_{2}])}
m([f^{i}(y_{1}),f^{i}(y_{2})])\\\\

   & \leq e^{2V}\gamma \sum_{i\in I}m([f^{i}(y_{1}),f^{i}(y_{2})])\\
 & \leq e^{2V}\gamma \\
& \leq \eta_{0}.
 \end{array}
$ \\

Consequently

\begin{equation}\label{(1)}
    \sum_{i\in I} |\textrm{Df}(f^{i}(z_{2}))-\textrm{Df}(a_{i})|\leq \varepsilon
\end{equation}
\medskip

Also, $\big([f^{i}(z_{2}),b_{i}] \big)_{i\in I}$ are disjoint intervals of continuity of $\textrm{Df}$, so by Lemma \ref{l:31}, they satisfy:\\

$\begin{array}{ll}
   \sum_{i\in I}m([f^{i}(z_{2}),b_{i}]) & \leq \sum_{i\in I}m([f^{i}(z_{2}),f^{i}(z_{3})])
\\\\

   & \leq \sum_{i\in I} e^{2V} \dfrac{m([z_{2},z_{3}])}{m([y_{1},y_{2}])}
m([f^{i}(y_{1}),f^{i}(y_{2})])\\\\

   & \leq e^{2V}\gamma \sum_{i\in I}m([f^{i}(y_{1}),f^{i}(y_{2})])
\leq \eta_{0}.
 \end{array}
$ \\\\
Consequently

\begin{equation}\label{(2)}
    \sum_{i \in I} |\textrm{Df}(f^{i}(z_{2}))-\textrm{Df}(b_{i})| \leq \varepsilon
\end{equation}
\medskip

Combining (\ref{(1)}) and (\ref{(2)}), we get

\begin{equation}\label{(3)}
    \sum_{i \in I} |\textrm{Df}(b_{i})-\textrm{Df}(a_{i})|\leq 2\varepsilon .
\end{equation}

Set $D^{\ast}_{i}(f):= \max\big(D_{i}(f),(D_{i}(f))^{-1}\big)$ and
$P^{\ast}:= \prod_{i\in I} D^{\ast}_{i}(f)$.
\medskip

Since $\textrm{Df}$ takes its values in $[m_{1},M_{1}]~~(m_{1}>0)$, the following properties hold:\\

$0\leq D^{\ast}_{i}(f)-1 \leq \dfrac{1}{m_{1}}
|\textrm{Df}(b_{i})-\textrm{Df}(a_{i})|$.\\

$\sum_{i \in I} |D^{\ast}_{i}(f)-1|\leq \dfrac{2}{m_{1}}
\varepsilon$.\\

$ |\log D_{i}(f)|= \log D^{\ast}_{i}(f)\leq D^{\ast}_{i}(f)-1$\\

$\dfrac{1}{P^{\ast}} \leq P \leq P^{\ast}$.\\

It follows that
\begin{equation}\label{(4)}
    e^{-C_{0}\varepsilon} \leq P \leq e^{C_{0}\varepsilon}
 \, \, \ \textrm{where}  \, \, C_{0} = \dfrac{2}{m_{1}}.
\end{equation}
\medskip

\textbf{Step 2}. We consider $D_{k_{n}(d)}(f),~d\in E(c_{1},c)$. We have $d \in ~ ]f^{k_{n}(d)}(z_{1}),f^{k_{n}(d)}(z_{3})[$.\\

\textbf{Case 2a}. $d\in ]f^{k_{n}(d)}(z_{1}),f^{k_{n}(d)}(z_{2})[$
with $k_{n}(d)=i_{n}(d)$. By Proposition \ref{p:42} (5),
$\textrm{Df}$ is continuous on the intervals
$[f^{k_{n}(d)}(z_{1}),d],~[d,f^{k_{n}(d)}(z_{2})]$ and
$[f^{k_{n}(d)}(z_{2}),f^{k_{n}(d)}(z_{3})]$; by the mean value
theorem there exist $t_{1}\in ]f^{k_{n}(d)}(z_{1}),d[,~t_{2} \in
]d,f^{k_{n}(d)}(z_{2})[$ and $t_{3} \in
]f^{k_{n}(d)}(z_{2}),f^{k_{n}(d)}(z_{3})[$ such that
$$D_{k_{n}(d)}(f) = (1-r_{n}(d))\dfrac{\textrm{Df}(t_{1})}{\textrm{Df}(t_{3})}+r_{n}(d)\dfrac{\textrm{Df}(t_{2})}{\textrm{Df}(t_{3})}$$
~~where~~$$r_{n}(d) =
\dfrac{m([d,f^{k_{n}(d)}(z_{2})])}{m([f^{k_{n}(d)}(z_{1}),f^{k_{n}(d)}(z_{2})])}.$$

We have
\medskip

 $\begin{array}{ll}
           |D_{k_{n}(d)}(f)-\sigma_{f}(d)| & \leq |\dfrac{\textrm{Df}(t_{1})}{\textrm{Df}(t_{3})}-\dfrac{D_{-}f(d)}{D_{+}f(d)}|+
r_{n}(d)\dfrac{|\textrm{Df}(t_{1})-\textrm{Df}(t_{2})|}{\textrm{Df}(t_{3})}\end{array}
$
\medskip

Since $$\begin{array}{ll}
|\dfrac{\textrm{Df}(t_{1})}{\textrm{Df}(t_{3})}-\dfrac{D_{-}f(d)}{D_{+}f(d)}|
\leq \dfrac{M_{1}}{m_{1}^{2}} \big(
|\textrm{Df}(t_{1})-D_{-}f(d)|+|\textrm{Df}(t_{3})-D_{+}f(d)|\big),
\end{array}$$
\medskip

and by Lemma \ref{l:31}
\medskip

$$\begin{array}{ll}
 r_{n}(d)     & = \dfrac{m \big(f^{k_{n}(d)}([f^{-k_{n}(d)}(d),z_{2}]) \big)}{m \big(f^{k_{n}(d)}([z_{1},z_{2}]
\big)}\\\\
    & \leq e^{2V}\dfrac{m ([f^{-k_{n}(d)}(d),z_{2}]) }{m ([z_{1},z_{2}])} \\\\
& \leq e^{2V}\dfrac{m ([y_{1},y_{2}]) }{m ([z_{1},z_{2}])}\times
\dfrac{m ([f^{-k_{n}(d)}(d),z_{2}]) }{m ([y_{1},y_{2}])} \\\\
& \leq \dfrac{e^{2V}}{\beta} s_{n}(d)
 \end{array}$$

with ~~~~
$$s_{n}(d):= \dfrac{m ([f^{-k_{n}(d)}(d),z_{2}]) }{l
([y_{1},y_{2}])}.$$
\medskip

It follows that:

$$|D_{k_{n}(d)}(f)-\sigma_{f}(d)| \leq K_{0} \big(
\dfrac{s_{n}(d)}{\beta}+|\textrm{Df}(t_{1})-D_{-}f(d)|+|\textrm{Df}(t_{3})-D_{+}f(d)|
 \big)$$
where

$K_{0}=\max(\dfrac{M_{1}}{m_{1}}e^{2V},\dfrac{M_{1}}{m_{1}^{2}})$.
\medskip

Since $\underset{n\to +\infty}\lim s_{n}(d)=0$ (Corollary
\ref{c:37}), there exists $n_{d,\beta} \in \mathbb{N}_{0}$ such that
for any $n \in \mathbb{N}_{0},~ n \geq n_{d,\beta}$, one has
$s_{n}(d)< \beta\varepsilon $. The intervals
$(]t_{1},d[;~]d,t_{3}[)_{d\in E(c_{1},c)}$ are disjoint intervals of
continuity of $\textrm{Df}$,  they satisfy  by Lemma \ref{l:31}:
\bigskip

 $
\begin{array}{ll}
 m([t_{1},d])+m([d,t_{3}]) & \leq 2 m([f^{k_{n}(d)}(z_{1}),f^{k_{n}(d)}(z_{2})])+m([f^{k_{n}(d)}(z_{2}),f^{k_{n}(d)}(z_{3})])\\\\

 & \leq 3 e^{2V}\gamma m([f^{k_{n}(d)}(y_{1}),f^{k_{n}(d)}(y_{2})])\\\\

  & \leq 3 e^{2V} \gamma < \eta_{0}
\end{array}
$

 So, $$ |\textrm{Df}(t_{1})-D_{-}f(d)|+|\textrm{Df}(t_{3})-D_{+}f(d)|< \varepsilon $$
\medskip

Hence, there exists $n_{d,\beta} \in \mathbb{N}_{0}$ such that for
every $n \in \mathbb{N}_{0}$, $n \geq n_{d,\beta}$,

\begin{equation}\label{(5)}
    |D_{k_{n}(d)}(f)-\sigma_{f}(d)|\leq 2 K_{0}\varepsilon
\end{equation}
\smallskip

\textbf{Case 2b}. $d\in ]f^{k_{n}(d)}(z_{2}),f^{k_{n}(d)}(z_{3})[$
with $k_{n}(d)=j_{n}(d)$. By Proposition \ref{p:42} (b-5),
$\textrm{Df}$ is continuous on the intervals
$[f^{k_{n}(d)}(z_{1}),f^{k_{n}(d)}(z_{2})],~[f^{k_{n}(d)}(z_{2}),d]$
and $[d,f^{k_{n}(d)}(z_{3})]$. By the mean value theorem, there
exist $s_{1} \in ]d, f^{k_{n}(d)}(z_{3})[,~s_{2} \in
]f^{k_{n}(d)}(z_{2}),d[$ and $s_{3} \in
]f^{k_{n}(d)}(z_{1}),f^{k_{n}(d)}(z_{2})[$ such that
$$D_{k_{n}(d)}^{-1}(f)=(1-r^{\prime}_{n}(d))\dfrac{\textrm{Df}(s_{1})}{\textrm{Df}(s_{3})}+r^{\prime}_{n}(d)\dfrac{\textrm{Df}(s_{2})}{\textrm{Df}(s_{3})}$$
where
$$r^{\prime}_{n}(d) = \dfrac{m[f^{k_{n}(d)}(z_{2}),d])}{m[f^{k_{n}(d)}(z_{2}),f^{k_{n}(d)}(z_{3})])}.
$$
\medskip

As in Case 2a, we have
\medskip

 $\begin{array}{ll}
           |D_{k_{n}(d)}(f)^{-1}-\sigma_{f}(d)^{-1}| & \leq |\dfrac{\textrm{Df}(s_{1})}{\textrm{Df}(s_{3})}-\dfrac{D_{+}f(d)}{D_{-}f(d)}|+
r^{\prime}_{n}(d)
\dfrac{|\textrm{Df}(s_{1})-\textrm{Df}(s_{2})|}{\textrm{Df}(s_{3})}\end{array}
$
\bigskip

Since $\begin{array}{ll}
|\dfrac{\textrm{Df}(s_{1})}{\textrm{Df}(s_{3})}-\dfrac{D_{+}f(d)}{D_{-}f(d)}|\leq
\dfrac{M_{1}}{m_{1}^{2}}
\big(|\textrm{Df}(s_{1})-D_{+}f(d)|+|\textrm{Df}(s_{3})-D_{-}f(d)|
\big),\end{array} $ \
\\
\medskip

and by Lemma \ref{l:31}
\medskip

$$\begin{array}{ll}
 r^{\prime}_{n}(d) & = \dfrac{m \big(f^{k_{n}(d)}([z_{2},f^{-k_{n}(d)}(d)]) \big)}{m \big(f^{k_{n}(d)}([z_{2},z_{3}]
\big)}\\\\
    & \leq e^{2V}\dfrac{m ([z_{2},f^{-k_{n}(d)}(d)]) }{m ([z_{2},z_{3}])} \\\\
& \leq e^{2V}\dfrac{m ([y_{1},y_{2}]) }{m([z_{2},z_{3}])}\times
\dfrac{m ([z_{2},f^{-k_{n}(d)}(d)]) }{m ([y_{1},y_{2}])} \\\\
& \leq \dfrac{e^{2V}}{\beta} s_{n}^{\prime}(d),\end{array}
$$ where $$s_{n}^{\prime}(d):= \dfrac{m ([z_{2},f^{-k_{n}(d)}(d)]) }{m
([y_{1},y_{2}])}.$$ \
\\

It follows that $$|D_{k_{n}(d)}(f)^{-1}-\sigma_{f}(d)^{-1}| \leq
K_{0} \big(
\dfrac{s_{n}^{\prime}(d)}{\beta}+|\textrm{Df}(s_{1})-D_{+}f(d)|+|\textrm{Df}(s_{3})-D_{-}f(d)|
 \big)$$
where $$K_{0} =
\max(\dfrac{M_{1}}{m_{1}}e^{2V},\dfrac{M_{1}}{m_{1}^{2}}).$$ By
Corollary \ref{c:37} : $\underset{n\to +\infty}\lim
s_{n}^{\prime}(d)=0$. Hence, there exists $n_{d,\beta} \in
\mathbb{N}_{0}$ such that for every $n\in \mathbb{N}_{0}, \ n \geq
n_{d,\beta}$,~~ $s_{n}^{\prime}(d)< \beta\varepsilon $. The
intervals $(]d,s_{1}[;~]s_{3},d[)_{d\in E(c_{1},c)}$ are disjoint
intervals of
continuity of $\textrm{Df}$, they satisfy, by Lemma \ref{l:31}, \\
\medskip

 $\begin{array}{ll}
 m([d,s_{1}])+m([s_{3},d]) & \leq    m([f^{k_{n}(d)}(z_{1}),f^{k_{n}(d)}(z_{2})])+ 2m([f^{k_{n}(d)}(z_{2}),f^{k_{n}(d)}(z_{3})])\\\\

 & \leq 3e^{2V}\gamma m([f^{k_{n}(d)}(y_{1}),f^{k_{n}(d)}(y_{2})])\\\\

  & \leq 3e^{2V}\gamma < \eta_{0}
\end{array}
$

 So, $$|\textrm{Df}(s_{1})-D_{+}f(d)|+|\textrm{Df}(s_{3})-D_{-}f(d)|< \varepsilon $$
\medskip

Hence, there exists $n_{d,\beta}\in \mathbb{N}_{0}$ such that for
every $n \in \mathbb{N}_{0}, n\geq n_{d,\beta}$,

\begin{equation}\label{(6)}
    |D_{k_{n}(d)}(f)^{-1}-\sigma_{f}(d)^{-1}| \leq 2 K_{0}\varepsilon
\end{equation}
\medskip

Therefore from the cases 2a and 2b, we conclude that there exists
$n_{d,\beta} \in \mathbb{N}_{0}$ such that for every $n \in
\mathbb{N}_{0},~n \geq n_{d,\beta}$,
\medskip

\begin{equation}\label{(7)}
     |D_{k_{n}(d)}(f)-\sigma_{f}(d)| \leq  K_{0}\varepsilon
\end{equation}
\medskip

Since $E(c_{1},c)$ is finite, there exists $n_{\beta} \in
\mathbb{N}_{0}$ such that for every $n \in \mathbb{N}_{0}, ~n \geq
n_{\beta}$,

\begin{equation}\label{(8)}
   |\prod_{d \in E(c_{1},c)}~D_{k_{n}(d)}(f)-\prod_{d \in
E(c_{1},c)}~\sigma_{f}(d)| \leq \varepsilon
\end{equation}
\medskip

Hence, (\ref{(4)}) and (\ref{(8)}) imply that: there exist $n_{\beta} \in \mathbb{N}_{0}$ such that for every $n \in
\mathbb{N}_{0},~n \geq n_{\beta}$,

\begin{equation}\label{(9)}
    |Dcr_{f^{q_{n}}}(z_{1},z_{2},z_{3})-\nu(c) | \leq C_{1}\varepsilon
\end{equation}
\medskip

where $\nu(c) = \prod_{d\in E(c_{1},c)}~\sigma_{f}(d)$ and $C_{1}$ is a
positive constant.
\end{proof}
\bigskip

\begin{lem} 
\label{l:49} Under the hypothesis of Proposition \ref{p:44}, for any $\varepsilon >0$, there exists
$0<\gamma<\gamma_{0}$, such that for any $\beta \in ]0,\gamma[$,
there exist $n_{\beta}\in \mathbb{M}_{0}$ such that for any $n\in
\mathbb{M}_{0}, ~n\geq n_{\beta}$, the ($\beta,\gamma$)-secondary
cell $(z_{1},z_{2},z_{3})$ associated to $(f,x_{0},b,\mathbb{M}_{0}, \delta)$ satisfies

$$\arrowvert\textrm{Dr}_{g^{q_{n}(b)}}(h(z_{1}),h(z_{2}),h(z_{3}))-
\Pi_{g}(h(c), h(c_{1}))\arrowvert\leq C_{2}\varepsilon$$
\medskip

where \ $\Pi_{g}(h(c), h(c_{1})) = \underset{d\in E (c_{1},c)}\prod \sigma_{g}(h(d))$ and \; $C_{2}$ is a
positive constant.
\end{lem}
\medskip

\begin{proof} The proof is a consequence of Lemma \ref{l:49q} and (Proposition \ref{p:43}, (h-5)) applied to
the ($\frac{\beta}{2},\gamma$)-derived cell
$(h(z_{1}),h(z_{2}),h(z_{3}))$ associated to $(g,h(x_{0}),h(c),
\mathbb{N}_{0})$ instead of the ($\beta,\gamma$)-derived cell
$(z_{1},z_{2},z_{3})$ associated to $(f,x_{0},c,\mathbb{N}_{0})$.
 \end{proof}
\medskip
\
\\
{\it Proof of Proposition \ref{p:44}}. The proof results from the
Lemmas \ref{l:49q} and \ref{l:49}. \qed
\bigskip

   \section{\bf Proof of Main Theorem }
   \medskip

Let $f, ~g\in \mathcal{P}(S^{1})$ with the same irrational rotation number $\alpha$ of bounded type.
Suppose that $f$ and $g$ satisfy the (KO) condition. By Corollary \ref{c:21}, there exist two
 piecewise quadratic homeomorphisms $K, L\in \mathcal{P}(S^{1})$ such that $F= L \circ f \circ L^{-1}$ and
 $G = K \circ g \circ K^{-1}$ have the following properties:

   \begin{enumerate}
     \item $F, G\in \mathcal{P}(S^{1})$ and have the same irrational rotation number $\alpha$,

     \item  The break points of $F$ (resp. $G$) are on \textit{pairwise distinct} $F$-orbits (resp. $G$-orbits),

     \item $F$ and $G$ satisfy the (KO) condition,

     \item Let $h$ the conjugating map between $f$ and $g$ ~ i.e. $h \circ f = g \circ h$ and set $v = K \circ h \circ L^{-1}$. Then

     \subitem - $v$ is an absolutely continuous (resp. a singular) function if and only if so is $h$.

    \subitem - $v\circ F= G\circ v$.
   \end{enumerate}

   Therefore, we may assume that all break points of $f$ (resp. $g$) are on \textit{pairwise distinct} $f$-orbits
   (resp. $g$-orbits). Now by Proposition \ref{p:23}, there exists a
   homeomorphism $u$ of $S^{1}$ such that $G= u\circ f\circ u^{-1}\in  \mathcal{P}(S^{1})$ with the following properties:

   \
   \\
   \textbullet ~ $C(f)\subset B:= \{c_{i}: ~i=0,1, \dots, q\}$, ($q\geq p$). \\
   \textbullet ~ $C(G)\subset u(B):=\{u(c_{i}): ~i=0,1, \dots, q\}$. \\
   \textbullet ~ $\sigma_{G}(u(c_{i})) = \sigma_{g}(d_{i})~,~i=0,1, \dots, q$.\\
   \textbullet ~ $\pi_{s}(G)=\pi_{s}(g)$.\\
   \textbullet ~ $h$ is singular if and only if so is $u$.
   \medskip

   So we may assume that $u=h$ and $G=g$.
\
\\
   To show that the conjugation homeomorphism $h$ from $f$ to $g$ is singular with respect to the Lebesgue measure $m$, it
   suffices to prove that its derivative $\textrm{Dh}$ is zero on
   a set of Lebesgue total measure. Assume on the contrary that $h$ admits at a point $x_{0}$ a positive derivative $Dh(x_{0})>0$.
\bigskip

Let $\varepsilon >0, ~c \in B$. For $0<\beta<\gamma <\gamma_{0}$,
write:

$$D_{n}(\beta, \gamma):= \dfrac{\textrm{Dr}_{g^{q_{n}}}\big(h(z_{1}),h(z_{2}), h(z_{3}) \big)}{\textrm{Dr}_{f^{q_{n}}}\big(z_{1},z_{2}, z_{3} \big)} $$
\medskip

By Proposition \ref{p:2}, there exists $n_{\beta} \in
\mathbb{N}_{0}$ such that for every $n\in \mathbb{N}_{0},~n\geq
n_{\beta}$:

\begin{equation}\label{(10)}
    |D_{n}(\beta, \gamma)-1| \leq \dfrac{2}{1-\gamma_{0}}\beta
\end{equation}
\medskip

By Proposition \ref{p:44}, there exists
$\gamma_{\varepsilon}<\gamma_{0}$ such that for every
$0<\beta<\gamma_{\varepsilon}$ there exists
$n_{\beta}(\varepsilon)\in \mathbb{N}_{0}$ such that for every $n\in
\mathbb{N}_{0},~~n\geq n_{\beta}(\varepsilon)$, the $(\beta,\gamma)$
secondary cell $(z_{1},z_{2},z_{3})\ (\gamma<\gamma_{\varepsilon})$-
associated to $(f,x_{0},c,\mathbb{N}_{0},\delta)$ satisfies:

\begin{equation}\label{(11)}
   |D_{n}(\beta, \gamma)-\Pi(c_{1},c)| \leq A \varepsilon,
\end{equation}
\
\\
where $$\Pi(c_{1},c):= \underset{d \in
E(c_{1},c)}\prod~\dfrac{\sigma_{g}(h(d))}{\sigma_{f}(d)}$$ and $A$
is a positive
constant.\\

Let $\beta < \gamma <\gamma_{\varepsilon}$ and $n \in
\mathbb{N}_{0}, ~n \geq n_{\beta}(\varepsilon)$. Then (\ref{(10)})
and (\ref{(11)}) imply that
\bigskip

$ \begin{array}{lll}
            |\Pi(c_{1},c)-1| & \leq |D_{n}(\beta, \gamma)-1|+|D_{n}(\beta,
\gamma)-\Pi(c_{1},c)| & \leq \dfrac{2\beta}{1-\gamma_{0}}+A
\varepsilon

  \end{array}$
\bigskip

Since $\varepsilon$ and $ \beta $ are arbitrary, it follows that $\Pi(c_{1},c)=1$. i.e.
\begin{equation}\label{(12)}
   \prod_{d\in E(c_{1},c)} ~\dfrac{\sigma_{g}(h(d))}{\sigma_{f}(d)}=1
\end{equation}
Since $c_{1}, c \in B$,~~$c_{1}\neq c~$ are arbitrary, so for every
$(i,k)\in \{0,1, \dots,~q \}^{2}~~(i \neq k),$ one has
$$\underset{d\in E(c_{k}, c_{i})}\prod~\dfrac{\sigma_{g}(h(d))}{\sigma_{f}(d)}=1$$

As by Proposition \ref{p:461}, $c_{k} \notin E(c_{k},c_{i})$, so we have $$\underset{0 \leq j \leq q
}\prod~\big (\dfrac{\sigma_{g}(h(c_{j}))}{\sigma_{f}(c_{j})}
\big)^{e(i,k,j)}=1,$$

where $e(i,k,j) \in \{0,1\}$ with
$$\varepsilon(i,k,i)=1 \textrm{and } \varepsilon(i,k,k)=0.$$
It follows that for every $(i,k)\in \{0,1, \dots,~q \}^{2}, \ i\neq
k$, we have

$$\underset{0 \leq j \leq q}\sum~ e (i,k,j)~\log\left(\dfrac{\sigma_{g}(h(c_{j}))}{\sigma_{f}(c_{j})}\right)=0.$$
\medskip
 \
\\

For $(i,j,k) \in \{0, \dots, q\}^{3}$, set
$$\varepsilon_{k+i(q+1),j}
    =\left\{
                            \begin{array}{ll}
                         e (i,k,j) & \textrm{if} ~~i \neq k,~~j \neq i ~~\textrm{and}~~j\neq
k,
\\\\
                       1 &  \textrm{if} ~~i \neq k ~~\textrm{and}~~j=
i,
\\\\
                          0 & \textrm{if} ~~i \neq k
~~\textrm{and}~~j=k,
\\\\
0   &   \textrm{if} ~~i = k
                   \end{array}
         \right.
$$

It follows that for every $(i,k) \in \{0, \dots,
q\}^{2}$,~~$$\sum_{j=0}^{q}\varepsilon_{k+i(q+1),j}\log\left(\frac{\sigma_{g}(h(c_{j}))}{\sigma_{f}(c_{j})}\right)=0
$$

Set $A_{q+1}=(a_{l,j})_{0\leq l<(q+1)^{2},~0 \leq j<q}\in
\mathcal{M}_{(q+1)^{2},~q+1}(\{0,1\})$  the matrix defined by
$a_{l,j}=\varepsilon_{k+i(q+1),j}~$, where $l=k+i(q+1)$ with $0\leq
i,~k \leq q$. The matrix $A_{q+1}$ has a rank $q+1$ (see \cite{am12}, Lemma 6.1). Therefore for every
$j\in \{0,1, \dots,~q\}$, $\sigma_{g}(h(c_{j}))=\sigma_{f}(c_{j})$. Hence $f$ and $g$ are break-equivalent, a contradiction. 
The proof is complete. \qed
\
\\

\textit{Proof of Corollary \ref{c:12}}. This follows from the fact that if
$f$ and $g$ do not have the same number of singular orbits then $f$ and
$g$ are not break-equivalent. \qed
\
\\

\textit{Proof of Corollary  \ref{c:13}}. Suppose that $f$ and $g$ are
break-equivalent. After reduction, one can suppose that $\textrm{SO}(f)= C(f)$ and $\textrm{SO}(g)= C(g)$. Then for every $d \in C(g)$, there is $e \in C(f)$
such that $d=h(e)$ and $\sigma_{g}(d)=\sigma_{f}(e)$ i.e.
$\sigma_{g}(d) \in \{\sigma_{f}(c):~c \in C(f)\}$, this
contradicts the hypothesis of the corollary.\qed
\
\\

\textit{Proof of Corollary  \ref{c:14}}. This follows from the fact that the
condition (ii) satisfies the assumption of Corollary \ref{c:13}. \qed
\
\\

\textit{Proof of Corollary  \ref{c:15}}. \textit{Proof of (i)}. Since $g$
does not have the $(D)-$ property and $f$ has the $(D)-$ property,
so there exists a point $d \in C(g)$ such that
$\pi_{s,O_{g}(d)}(d)\neq 1$ and $\pi_{s,O_{f}(c)}(f)= 1$ for every
$c \in C(f)$. We conclude by Corollary \ref{c:13}.
\
\\

 \textit{Proof of (ii)}. By Corollary \ref{c:21}, $f$ is conjugate
to a class P-homeomorphism $F=K \circ f \circ K^{-1}$ through a
piecewise quadratic class $P$ homeomorphism $K$ of $S^{1}$. Since
$f$ satisfies the $D-$ property, $\pi_{s,O_{f}(c)}(f)=
1=\sigma_{F}(K(c))$ for every $c \in C(f)$. Thus, $F$ is a $C^{1}$-
diffeomorphism. Moreover, $\textrm{DF}$ is absolutely continuous on $S^{1}$
and $D^{2}F \in L^{p}(S^{1})$ for some $p>1$. By
Katznelson-Ornstein's theorem, $\mu_{F}$ is equivalent to the
Lebesgue measure $m$ and so is $\mu_{f}$. Hence there exists an
absolutely continuous map $\varphi$ such that $\varphi \circ f=
R_{\rho(f)} \circ \varphi$. Similarly, there exists an absolutely
continuous map $\psi$ such that $\psi \circ g = R_{\rho(f)} \circ
\psi$. In addition, $\psi^{-1}$ is absolutely continuous. It
follows that $(\psi \circ h\circ\varphi^{-1}) \circ
R_{\rho(f)}=R_{\rho(f)}\circ (\psi\circ h\circ \varphi^{-1})$. As
$\rho(f)$ is irational then $\psi\circ h\circ
\varphi^{-1}=R_{\beta}$ a rotation, for some $\beta \in S^{1}$.
Therefore $h= \psi^{-1}\circ R_{\beta}\circ \varphi$, which is an
absolutely continuous map. This completes the proof. \qed
\
\\

\textit{Proof of Corollary  \ref{c:16}}. The proof follows easily from
Corollary \ref{c:15} by taking $g$ the rotation $R_{\rho(f)}$. \qed



\bigskip

\bibliographystyle{amsplain}

\end{document}